# The Cloven Traveling Salesman: Cycle Covers and the Integrality Gap of Small ATSP Instances*

Alessandro Sosso[1], Ambrogio Maria Bernardelli[2,†], and Stefano Gualandi[2]

[1]Department of Computer Science, Aarhus University, Åbogade 34, Aarhus N, 8200, Denmark

[2]Department of Mathematics "F. Casorati", University of Pavia, via Ferrata 5, 27100 Pavia, Italy

[†]Corresponding author: ambrogiomaria.bernardelli@unipv.it

**Abstract**

This work proposes a novel enumeration algorithm for computing the integrality gap of small instances of the subtour elimination formulation for the Asymmetric Traveling Salesman Problem (ATSP). The core idea is to enumerate pairs of cycle covers that can be filtered and mapped to half-integer vertices of the subtour elimination polytope. The two-cycle covers are encoded as lexicographically ordered partitions of $n$ numbers, with an encoding that prevents the generation of several isomorphic vertices. However, since not every cycle cover pair can be mapped to a vertex of the subtour elimination polytope, we have designed an efficient property-checking procedure to control whether a given point is a vertex of the asymmetric subtour elimination polytope. The proposed approach turns upside down the algorithms presented in the literature that first generate every possible vertex and later filter isomorphic vertices. With our approach, we can replicate state-of-the-art results for $n \leq 9$ in a tiny fraction of time, and we compute for the first time the exact integrality gap of half-integer vertices of the asymmetric subtour elimination polytope for $n = 10, 11, 12$.

# 1 Introduction

The Traveling Salesman Problem (TSP) is a long-established combinatorial optimization problem consisting of finding the shortest tour visiting every node in a given edge-weighted graph [12]. Given $n$ nodes and arc costs $c_{ij}$ of traveling from node $i$ to node $j$, we focus on *metric* costs, which are costs that obey the triangle inequality $c_{ik} \leq c_{ij} + c_{jk}$, for all $i, j, k$. If the costs are also symmetric, we have the Symmetric TSP (STSP); otherwise, if we do not require $c_{ij} = c_{ji}$, we deal with the Asymmetric TSP (ATSP).

*This version of the article has been accepted for publication, after peer review, but is not the Version of Record and does not reflect post-acceptance improvements, or any corrections. The Version of Record is available online at: http://dx.doi.org/10.1007/s11590-026-02312-x.

Among the several algorithms designed over time to solve the TSP, a key method was developed by Dantzig, Fulkerson, and Johnson [13], who introduced a fundamental formulation of the TSP as an Integer Linear Problem (ILP). Despite considerable achievements in improving the TSP algorithms, the inherent difficulty of finding an efficient polynomial time algorithm led many efforts to be devoted to finding algorithms that could provide approximate solutions in polynomial time, both for the symmetric [10, 22] and for the asymmetric case [1, 16, 26, 27, 28]. From an integer programming perspective, a natural source of performance guarantees is the linear relaxation of the TSP formulation, obtained by dropping the integrality constraints; we refer to this relaxation as the Subtour Elimination Problem (SEP). For minimization problems, the optimal value of the linear relaxation provides a lower bound on the optimal integer value. A standard measure of the quality of this bound is the integrality gap, defined as the ratio between the optimal value of the integer problem and that of its linear relaxation (e.g., see [11]). An upper bound on the integrality gap, together with a polynomial-time procedure that converts a relaxation solution into a feasible tour with the same worst-case guarantee, yields an approximation algorithm with the same ratio.

While the integrality gap of the TSP is known to be unbounded in the general case, since instances achieving any arbitrary gap can be constructed [15, S4.2], the computation of the maximum integrality gap when considering only metric weights is an open question, as only lower and upper bounds have been proved. Structural analyses of the subtour LP, such as those by Carr and Vempala [7], have identified families of instances where the gap can be tightly characterized, helping to formalize how the LP behaves across different classes of graphs. For the metric STSP, the integrality gap of the subtour LP has long been known to be at most $\frac{3}{2}$ (e.g., see [31, 23]), and no instances with a gap of $\frac{4}{3}$ or greater have been found. Recently, Karlin et al. [20] showed that the integrality gap is strictly less than $\frac{3}{2}$, giving the first improvement on this bound in decades. In the asymmetric case, only a large upper bound of 22 was known [28], which has recently been improved to 17 [29], while the best lower bound has been discovered by Charikar, Goemans, and Karloff [9] and independently by Elliott-Magwood [15] with a family of weighted graphs with an integrality gap tending to 2, and analogous constructions for unweighted graphs are discussed in [29].

The computation of the integrality gap for small TSP instances (i.e., with less than 15 nodes), is of interest since it might provide useful insights into the structure of the instances yielding a high integrality gap. As the maximal integrality gap cannot be computed directly since it would require computing the gap for infinite instances of the TSP, indirect approaches are to be employed. One general framework for analyzing integrality gaps is the fractional decomposition tree method introduced by Carr et al. [8], which provides theoretical bounds by decomposing fractional solutions. In this work, however, we rely on a more computationally driven strategy, as exhibited in [2, 6, 15, 30], where an auxiliary LP problem, derived from the dual of the SEP, is leveraged. Such auxiliary GAP problem takes as input a feasible solution to the SEP, and by calculating its optimum on all vertices of the SEP polytope for a given $n$, the maximal value $\text{Gap}_n$ of the integrality gap over all $n$-nodes TSP instances can be obtained. This approach is not exclusive to the TSP problem, as a similar methodology was recently employed to compute the integrality gap of small metric Steiner Tree problems [3], thereby highlighting the generalizability of the approach.

The goal of this work is to revisit and extend the approach presented in [15] to compute for the first time the exact $\text{Gap}_n$ for $n = 10, 11, 12$, for all *half-integer* extreme points of ATSP, that are those vertices with component values in $\{0, 1/2, 1\}$. This choice is motivated by computational evidence, as this type of vertices achieve the highest known integrality gaps for those $n$ for which the exact $\text{Gap}_n$ has been computed over all extreme points [15, S4]. In particular, when talking about the integrality gap of an extreme point, we refer to the maximum ratio between the optimal value of the integer problem and the optimal value of the linear relaxation, with the additional constraint that the extreme point under consideration is an optimal solution of the relaxation. Notice that in [15] and in the previous work [5], the authors could compute $\text{Gap}_n$ for $n \leq 7$ for any type of vertex via exhaustive generation of all the vertices of ASEP. For $n = 8, 9$, however, the number of vertices was too large for their extensive computation, and the authors focused on computing the integrality gap of all *half-integer* extreme points. For $n > 9$, no form of complete enumeration was computationally practicable, but the authors nonetheless provided lower bounds on the value of $\text{Gap}_n$ by concentrating on only a few specific instances. To the best of our knowledge, these constitute the most advanced results to date on the exact computation of the maximum ATSP integrality gap for half-integer vertices.

This paper builds upon the master's thesis [24], where parts of the structural and algorithmic development appear in more detailed form. The present version refines the methodology and extends the computational results to $n = 12$.

*Our contributions.* This paper develops a novel enumeration framework for the exhaustive computation of half-integer extreme points of ATSP and their integrality gap. Our starting point is a structural observation: every half-integer extreme point can be represented as the average of two edge-disjoint cycle covers (Theorem 2). Here and throughout the paper, by a cycle cover we mean a spanning directed 1-regular subgraph, and therefore a union of vertex-disjoint directed cycles. Decomposition results for fractional solutions follow from known properties of circulation polytopes (see, e.g., [17, 29]); however, to the best of our knowledge this representation has not been formulated explicitly in this vertex-level and algorithmic form. We include a self-contained derivation in Section 2, as it forms the structural foundation of our enumeration procedure. Building on this structure, we design a new generation-and-filtering framework that constructs a superset of candidate vertices via pairs of cycle covers, encodes them to remove redundancies, and certifies them through feasibility and extremality checks. Unlike [15], where all vertices are first enumerated and only later analyzed, our approach interleaves structural generation and certification, significantly reducing the effective search space and enabling exact computations for larger instances. A key component of the pipeline is a novel circuit-based extremality verification procedure (Propositions 1 and 2) tailored to the combinatorial structure of the generated candidates, which substantially accelerates the extremality test. Combining these ingredients, we compute for the first time the integrality gap of all half-integer extreme points for $n \leq 12$, extending previous exact computations and confirming the conjectured maximum gap within this class. Compared to earlier exhaustive approaches, our framework achieves these results within a fraction of the previously required computation time. Beyond these specific computations, the structural generation and certification techniques introduced here provide new methodological perspectives on studying extreme-point structure, which may be applicable to related polyhedral and combinatorial optimization problems.

*Outline.* The outline of this paper is as follows. Section 2 presents the problem of computing the integrality gap of the ATSP and derives the structural characterization of half-integer vertices that underpins our enumeration algorithm. Section 3 and Section 4 detail exhaustive and efficient algorithms for generating all half-integer vertices for a fixed $n$. In particular, while Section 3 exploits Theorem 2 for generating candidate points that could represent half-integer vertices, Section 4 presents an efficient procedure for checking whether a given candidate point is a basic feasible solution, leveraging the structure of the problem and designing a novel algorithm. In Section 5, we report our computational results on the complete enumeration of half-integer vertices for $n \leq 12$, and the evaluation of the integrality gap of each vertex. We conclude the paper with a perspective on future work.

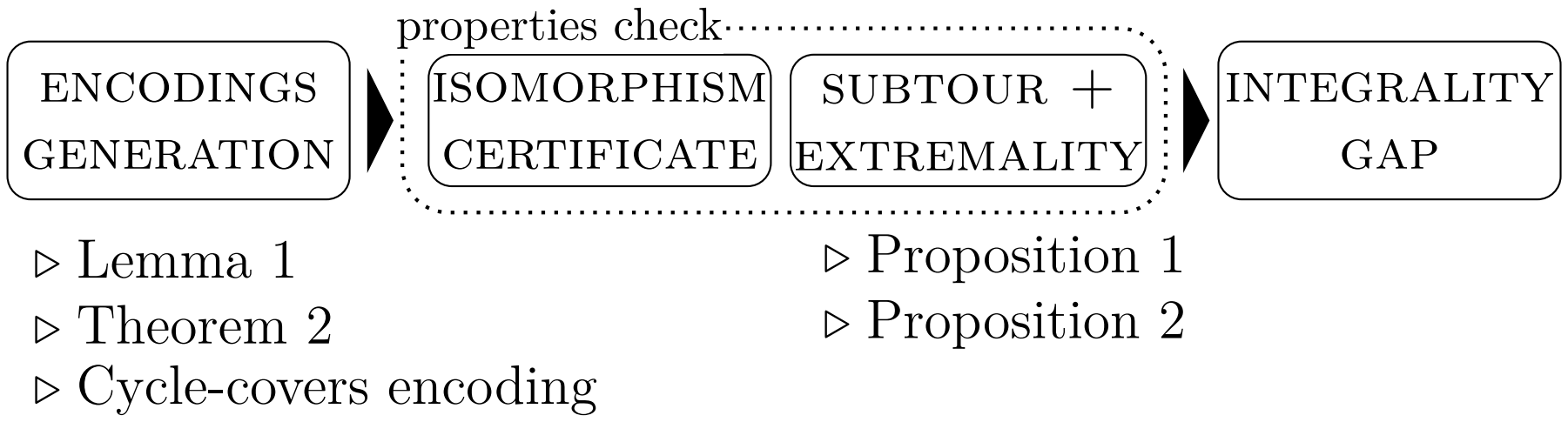

Figure 1: Overview of the proposed approach for computing the maximum integrality gap of half-integer vertices, with references to the results presented in this paper.

# 2 Gap problem and vertex characterization

## 2.1 The gap problem

Given a directed weighted graph $G = (V, E, \boldsymbol{c})$ with $n = |V|$, the classical mixed integer linear programming formulation of the ATSP uses 0–1 arc variables $x_{uv}$, with $(u, v) \in E$, to represent the characteristic vector of tours as follows:

$$
\begin{aligned}
(\text{ATSP}) \quad \min \quad & \textstyle\sum_{uv \in E} c_{uv} x_{uv} && && \text{(1a)}\\
\text{s.t.} \quad & \textstyle\sum_{uv \in \delta^+(w)} x_{uv} = 1 && \forall w \in V, && \text{(1b)}\\
& \textstyle\sum_{uv \in \delta^-(w)} x_{uv} = 1 && \forall w \in V, && \text{(1c)}\\
& \textstyle\sum_{uv \in \delta^+(S)} x_{uv} \geq 1 && \forall S \in \mathcal{S}^n, && \text{(1d)}\\
& x_{uv} \geq 0 && \forall uv \in E, && \text{(1e)}\\
& x_{uv} \in \mathbb{Z} && \forall uv \in E, && \text{(1f)}
\end{aligned}
$$

with $\mathcal{S}^n := \{S \subset V \mid 2 \leq |S| \leq n-2\}$. The resulting LP relaxation of ATSP is the Asymmetric Subtour Elimination Problem (ASEP). Note that any instance and solution

of both ATSP and ASEP is determined by the vector of costs $\boldsymbol{c}$, therefore we write the optimal value as ATSP $(\boldsymbol{c})$, ASEP $(\boldsymbol{c})$, respectively. Hence, we denote by $\mathfrak{P}^n$ the polytope that represents the region of the feasible solutions to the ASEP for a (complete) digraph on $n$ vertices $\mathcal{K}_n = (V^n, E^n)$, that is $\mathfrak{P}^n := \{\boldsymbol{x} \in \mathbb{R}^m \mid (1b), (1c), (1d), (1e)\}$, where $m := |E^n| = n(n-1)$. We also denote the (finite) set of extreme points of $\mathfrak{P}^n$ by $\mathfrak{X}^n$.

For the ATSP with a non-negative cost $\boldsymbol{c} \neq \boldsymbol{0}$, the integrality gap of the ASEP relaxation of ATSP is defined as $\text{Gap}\,(\boldsymbol{c}) := \text{ATSP}\,(\boldsymbol{c}) \,/\, \text{ASEP}\,(\boldsymbol{c}) \geq 1$. Since the integrality gap is unbounded in the general case, we restrict our discussion to metric ATSP weights. We introduce the following quantity, expressing the maximal value of the integrality gap over all $n$-nodes, metric, non-degenerate ATSP instances as $\text{Gap}_n := \max\left\{\text{Gap}\,(\boldsymbol{c}) \;\middle|\; \boldsymbol{c} \in \mathbb{R}^{n(n-1)},\ \boldsymbol{c} \text{ metric},\ \boldsymbol{c} \neq \boldsymbol{0}\right\}$.

To compute $\text{Gap}_n$, the computational approach proposed in [2, 6, 15, 30] is followed, which relies on the solution of an auxiliary LP problem herein called GAP. For a given extreme point $\overline{\boldsymbol{x}} \in \mathfrak{X}^n$ of $\mathfrak{P}^n$ and for $\mathcal{S}^n_{uv} := \left\{S \in \mathcal{S}^n \;\middle|\; uv \in \delta^+(S)\right\}$, the problem GAP, which can be derived from the dual of ATSP, is defined by

$$
\begin{array}{llll}
(\text{GAP}) & \min & \displaystyle\sum_{uv \in E^n} \overline{x}_{uv} c_{uv}, & \\
& \text{s.t.} & c_{uw} + c_{wv} - c_{uv} \geq 0 & \forall uv \in E^n, w \in V^n \setminus \{u, v\}, \\
& & \text{cost}\,(T) \geq 1 & \forall T \text{ tour in } \mathcal{K}_n, \\
& & c_{uv} - y^{\text{out}}_u - y^{\text{in}}_v - \sum_{S \in \mathcal{S}^n_{uv}} d_S \geq 0 & \forall uv \in E^n \setminus \text{supp}_{E^n}(\overline{\boldsymbol{x}}), \\
& & c_{uv} - y^{\text{out}}_u - y^{\text{in}}_v - \sum_{S \in \mathcal{S}^n_{uv}} d_S = 0 & \forall uv \in \text{supp}_{E^n}(\overline{\boldsymbol{x}}), \\
& & c_{uv} \geq 0 & \forall uv \in E^n, \\
& & d_S \geq 0 & \forall S \in \mathcal{S}^n.
\end{array}
$$

Following [5], this allows us to exploit the relation $1/\text{Gap}_n = \min_{\overline{\boldsymbol{x}} \in \mathfrak{X}^n} \left\{ \text{GAP}\,(\overline{\boldsymbol{x}}) \right\}$. In particular, $1/\,\text{GAP}\,(\overline{\boldsymbol{x}})$ is the integrality gap of vertex $\overline{\boldsymbol{x}}$, that is, the highest ratio between the optimal value of the integer problem and the optimal value of the linear relaxation, with the additional constraint that $\overline{\boldsymbol{x}}$ is an optimal solution of the linear relaxation.

Even though this formulation makes computing the integrality gap possible, for relatively small values of $n$ (e.g., $8 \leq n \leq 12$) its actual solution becomes computationally intractable due to the exponential number of variables and constraints in GAP and for the number of points in $\mathfrak{X}^n$. Hence, we need an efficient approach to reduce the number of candidate points from $\mathfrak{X}^n$ to consider, enabling us to solve GAP for larger values of $n$.

## 2.2 Structure of half-integer solutions

Every point in $\mathfrak{X}^n$ is the solution to a system of linear equations with integer coefficients. Hence, any $\overline{\boldsymbol{x}} \in \mathfrak{X}^n$ has only rational components, that is $\mathfrak{X}^n \subset \mathbb{Q}^m \cap [0,1]^m$. This allows us to write any $\overline{\boldsymbol{x}} \in \mathfrak{X}^n$ as $\overline{\boldsymbol{x}} = (k_1/\alpha, k_2/\alpha, \ldots, k_m/\alpha)$ for some appropriate $\alpha \in \mathbb{N}^+, k_i \in \mathbb{N}, 0 \leq k_i \leq \alpha, i = \{1, \ldots, m\}$. Let us define

$$
\mathfrak{P}^n_\alpha := \mathfrak{P}^n \cap \{0, 1/\alpha, 2/\alpha, \ldots, 1\}^m, \qquad \mathfrak{X}^n_\alpha := \mathfrak{X}^n \cap \{0, 1/\alpha, 2/\alpha, \ldots, 1\}^m = \mathfrak{X}^n \cap \mathfrak{P}^n_\alpha,
$$

for which $\mathfrak{X}^n = \bigcup_{\alpha=1}^{\infty} \mathfrak{X}^n_\alpha$ holds. Herein, we focus on the case $\alpha = 2$, corresponding to half-integer vertices, as motivated in the introduction.

**Definition 1** (Half-integer and pure half-integer extreme points)**.** Let $\overline{\boldsymbol{x}} \in \mathfrak{X}^n$. We say that $\overline{\boldsymbol{x}} \in \mathfrak{X}^n$ is a *half-integer* extreme point if $\overline{\boldsymbol{x}} \in \mathfrak{X}_2^n$, that is if the components of $\overline{\boldsymbol{x}}$ take values in the set $\{0, 1/2, 1\}$. We say that $\overline{\boldsymbol{x}}$ is a *pure half-integer* extreme point if $\overline{\boldsymbol{x}} \in \mathfrak{X}_2^n \cap \{0, 1/2\}^m$.

A result by Elliott-Magwood [15, S4.6] proves that to compute the gap over all half-integer vertices it is enough to consider only pure half-integers, which greatly reduces the number of instances that need to be generated.

In the remainder of this section, we characterize half-integer solutions through two theoretical results, linking the points of $\mathfrak{X}_2^n$ to cycle covers. Let us first focus on the vectors associated with cycle covers in the complete graph $\mathcal{K}_n$, and let us call $\mathfrak{C}^n \subset \{0,1\}^m$ the set of these vectors. This vector space is characterized by the same constraints as the ATSP without the subtour elimination constraints (1d), $\mathfrak{C}^n = \{\overline{\boldsymbol{y}} \in \mathbb{R}^m \mid (1b), (1c), (1e), (1f)\}$. By additionally dropping the integrality constraints (1f), the set of extreme points of the resulting polytope $\mathfrak{Q}^n$ is equivalent to $\mathfrak{C}^n$, since Birkhoff showed in [4] how such extreme points are always integral. In Theorem 2, we claim that for every $\overline{\boldsymbol{x}} \in \mathfrak{X}_2^n$ two edge-disjoint cycle covers $\overline{\boldsymbol{y}}^1, \overline{\boldsymbol{y}}^2 \in \mathfrak{C}^n$ exist such that $\overline{\boldsymbol{x}} = \frac{1}{2}\overline{\boldsymbol{y}}^1 + \frac{1}{2}\overline{\boldsymbol{y}}^2$. To prove it, we first need to introduce an auxiliary lemma.

**Lemma 1.** *Let $\overline{\boldsymbol{x}} \in \mathfrak{X}_2^n$ and let $\overline{\boldsymbol{y}}^1 \in \mathfrak{C}^n$ be the characteristic vector of a cycle cover of the graph associated with $\overline{\boldsymbol{x}}$. If $\overline{\boldsymbol{y}}^2 \coloneqq 2\overline{\boldsymbol{x}} - \overline{\boldsymbol{y}}^1$, then $\overline{\boldsymbol{y}}^2 \in \mathfrak{C}^n$.*

*Proof.* The fact that $\overline{\boldsymbol{y}}^1$ is a cover of the $\overline{\boldsymbol{x}}$-associated graph implies that $\text{supp}_{E^n}(\overline{\boldsymbol{y}}^1) \subseteq \text{supp}_{E^n}(\overline{\boldsymbol{x}})$, which in turn means that $\mathbf{0} \leq \overline{\boldsymbol{y}}^1 \leq 2\overline{\boldsymbol{x}}$ and $\overline{\boldsymbol{y}}^2 = 2\overline{\boldsymbol{x}} - \overline{\boldsymbol{y}}^1 \geq \mathbf{0}$. The boundary constraints (1e) are thus satisfied by $\overline{\boldsymbol{y}}^2$, which also satisfies the out-degree constraints (1b) since $\forall w \in V^n$

$$\sum_{uv \in \delta^+(w)} \overline{y}^2_{uv} = \sum_{uv \in \delta^+(w)} \left(2\overline{x}_{uv} - \overline{y}^1_{uv}\right) = 2 \sum_{uv \in \delta^+(w)} \overline{x}_{uv} - \sum_{uv \in \delta^+(w)} \overline{y}^1_{uv} = 2 \cdot 1 - 1 = 1.$$

Analogously for the in-degree constraints (1c). The integrality constraints (1f) are trivially satisfied since both $2\overline{\boldsymbol{x}}$ and $\overline{\boldsymbol{y}}^1$ are integral. □

The proof of Theorem 2 relies on a result by Hall [19] known in the literature as Hall's Marriage Theorem, which provides a sufficient and necessary condition for the existence of a perfect matching on a bipartite graph.

**Theorem 1** (Hall's Marriage Theorem)**.** *There exists a $V_1$-perfect matching on a bipartite graph $\mathcal{G} = (V, E)$ with $V_1, V_2 \subset V$ as bipartite sets if and only if $|S| \leq |N_S|$, $\forall S \in \mathcal{P}(V_1)$, where $N_S \coloneqq \bigcup_{s \in S} \delta(s)$ is the set of neighbors for any $S \subseteq V_1$.*

**Theorem 2** (Structure of half-integer vertices)**.** *For each half-integer vertex of ASEP $\overline{\boldsymbol{x}}$, there exist two edge-disjoint cycle covers $\overline{\boldsymbol{y}}^1, \overline{\boldsymbol{y}}^2$ such that $\overline{\boldsymbol{x}} = \frac{1}{2}\overline{\boldsymbol{y}}^1 + \frac{1}{2}\overline{\boldsymbol{y}}^2$.*

*Proof.* First, we have to prove that, for each $\overline{\boldsymbol{x}} \in \mathfrak{X}_2^n$, there exist $\overline{\boldsymbol{y}}^1, \overline{\boldsymbol{y}}^2 \in \mathfrak{C}^n$ such that $\overline{\boldsymbol{x}} = \frac{1}{2}\overline{\boldsymbol{y}}^1 + \frac{1}{2}\overline{\boldsymbol{y}}^2$. Thanks to Lemma 1, it suffices to show that a cycle cover of the graph associated with $\overline{\boldsymbol{x}}$ exists. This is equivalent to the existence of a perfect matching on the bipartite graph

$$\widetilde{\mathcal{X}} \coloneqq \left(V^n \times \{0,1\}, \widetilde{E}^n_{\overline{\boldsymbol{x}}}\right), \qquad \text{with} \quad \widetilde{E}^n_{\overline{\boldsymbol{x}}} \coloneqq \{(u,0)(v,1) \mid uv \in \text{supp}_{E^n}(\overline{\boldsymbol{x}})\}.$$

Hall's Marriage Theorem 1 grants the existence of a $(V^n \times \{0\})$-perfect matching if $|\overline{S}| \le |N_{\overline{S}}|$ for any $\overline{S} \in \mathcal{P}(V^n \times \{0\})$. Let $\overline{S} \coloneqq S \times \{0\}$ with $S \subseteq V^n$. We see that

$$\bigcup_{u\in\overline{S}} \left\{ uv \in \widetilde{E}^n_{\overline{\boldsymbol{x}}} \;\middle|\; v \in \delta^+(u) \right\} \subseteq \bigcup_{v'\in N_{\overline{S}}} \left\{ u'v' \in \widetilde{E}^n_{\overline{\boldsymbol{x}}} \;\middle|\; u' \in \delta^-(v') \right\},$$

from which it follows

$$|\overline{S}| = \sum_{u\in\overline{S}} 1 = \sum_{u\in\overline{S}} \Big( \sum_{v\in\delta^+(u)} \overline{x}_{uv} \Big) \le \sum_{v'\in N_{\overline{S}}} \Big( \sum_{u'\in\delta^-(v')} \overline{x}_{u'v'} \Big) = \sum_{v'\in N_{\overline{S}}} 1 = |N_{\overline{S}}|,$$

that is, $|\overline{S}| \le |N_{\overline{S}}|$, proving the existence of the sought $(V^n \times \{0\})$-perfect matching. But since $|V^n \times \{0\}| = |V^n \times \{1\}| = |V^n|$, any such $(V^n \times \{0\})$-perfect matching is also a perfect matching, concluding the proof. □

**Remark 1.** Lemma 1 and Theorem 2 do not exploit that $\overline{\boldsymbol{x}} \in \mathfrak{X}^n_2$ satisfies the subtour elimination constraints (1d). The results therefore also hold for any $\overline{\boldsymbol{x}} \in \mathfrak{Q}^n \cap \{0, 1/2, 1\}^m$, that is, the half-integer solutions of the remaining constraints (1b), (1c), and (1e). This also means that the results of Theorem 2 can be extended via induction to $\mathfrak{X}^n_\alpha$ for $\alpha \ge 3$, proving that the elements $\overline{\boldsymbol{x}} \in \mathfrak{X}^n_\alpha$ are expressible as $\overline{\boldsymbol{x}} = \frac{1}{\alpha}\sum_{i=1}^{\alpha} \overline{\boldsymbol{y}}^i$ with $\overline{\boldsymbol{y}}^i \in \mathfrak{C}^n$, $i \in \{1, ..., \alpha\}$. In fact, similarly to Lemma 1 we see that if a cycle cover $\overline{\boldsymbol{y}}^\alpha \in \mathfrak{C}^n$ of $\mathcal{X}^n_{\overline{\boldsymbol{x}}}$ exists then $\overline{\boldsymbol{x}}' \coloneqq \frac{1}{\alpha-1}(\alpha\overline{\boldsymbol{x}} - \overline{\boldsymbol{y}}^\alpha)$ belongs to $\mathfrak{Q}^n \cap \{0, \frac{1}{\alpha-1}, \dots, 1\}^m$ (but not to $\mathfrak{X}^n_{\alpha-1}$ in general, as $\overline{\boldsymbol{x}}'$ may not satisfy the subtour elimination constraints). Hence from the inductive hypothesis we would have $\overline{\boldsymbol{x}}' = \frac{1}{\alpha-1}\sum_{i=1}^{\alpha-1} \overline{\boldsymbol{y}}^i$ with $\overline{\boldsymbol{y}}^i \in \mathfrak{C}^n$ for $i \in \{1, \dots, \alpha-1\}$ and therefore

$$\overline{\boldsymbol{x}} = \frac{\alpha-1}{\alpha}\overline{\boldsymbol{x}} + \frac{1}{\alpha}\overline{\boldsymbol{y}}^\alpha = \frac{\alpha-1}{\alpha} \cdot \frac{1}{\alpha-1}\sum_{i=1}^{\alpha-1} \overline{\boldsymbol{y}}^i + \frac{1}{\alpha}\overline{\boldsymbol{y}}^\alpha = \sum_{i=1}^{\alpha} \overline{\boldsymbol{y}}^i.$$

Theorem 2 implies that the enumeration of half-integer vertices in $\mathfrak{X}^n_2$ as the sum of edge-disjoint cycle covers is indeed exhaustive. But it leaves open the following question:

*Does it hold that, for any choice of* $\overline{\boldsymbol{y}}^1, \overline{\boldsymbol{y}}^2 \in \mathfrak{C}^n$, $\overline{\boldsymbol{x}} = \frac{1}{2}\overline{\boldsymbol{y}}^1 + \frac{1}{2}\overline{\boldsymbol{y}}^2$ *is in* $\mathfrak{X}^n_2$*?*

This is not guaranteed in the general case, despite $\overline{\boldsymbol{x}}$ trivially being in $\{0, 1/2, 1\}^m$. There are two possible cases leading to $\overline{\boldsymbol{x}} \notin \mathfrak{X}^n_2$:

- ▷ $\overline{\boldsymbol{x}} \notin \mathfrak{P}^n$, because it might violate the subtour elimination constraints (1d) (a trivial example of this can be constructed by taking $\overline{\boldsymbol{y}}^1 \in \mathfrak{C}^n \setminus \mathfrak{P}^n$ and $\overline{\boldsymbol{y}}^2 = \overline{\boldsymbol{y}}^1$, for which $\overline{\boldsymbol{x}} = \frac{1}{2}\overline{\boldsymbol{y}}^1 + \frac{1}{2}\overline{\boldsymbol{y}}^2 = \overline{\boldsymbol{y}}^1 \notin \mathfrak{P}^n$);
- ▷ $\overline{\boldsymbol{x}} \in \mathfrak{P}^n$ yet $\overline{\boldsymbol{x}} \notin \mathfrak{X}^n$, because it is not granted that $\overline{\boldsymbol{x}}$ is an extreme point of $\mathfrak{P}^n$ (Figure 2 shows an example of the latter case).

In general, it is only guaranteed that (*i*) $\overline{\boldsymbol{x}} \in \{0, 1/2, 1\}^m$, and (*ii*) $\overline{\boldsymbol{x}}$ satisfies (1b), (1c) and (1e), by taking the weighted sum of the same constraints that $\overline{\boldsymbol{y}}^1, \overline{\boldsymbol{y}}^2 \in \mathfrak{C}^n$ satisfy. It then holds that $\overline{\boldsymbol{x}} \in \mathfrak{Q}^n \cap \{0, 1/2, 1\}^m$ and $\overline{\boldsymbol{x}} \in \mathfrak{X}^n_2$ if and only if it satisfies the subtour elimination constraints (1d) and is an extreme point of $\mathfrak{P}^n$.

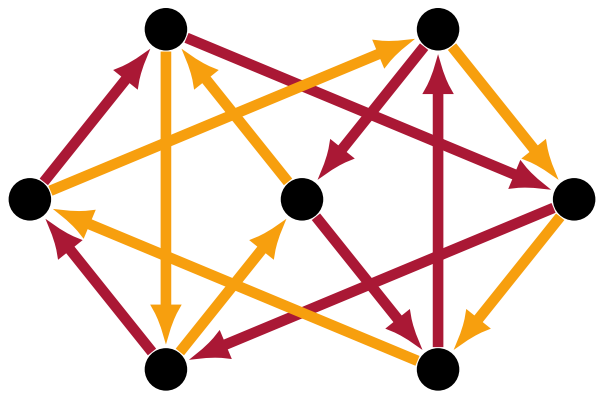

Figure 2: A solution $\overline{\boldsymbol{x}} = \frac{1}{2}\overline{\boldsymbol{y}}^1 + \frac{1}{2}\overline{\boldsymbol{y}}^2$ that is not an extreme point of $\mathfrak{P}^n$, since for all $S \in \mathcal{S}^n$ we have $\sum_{uv \in \delta^+(S)} x_{uv} > 1$, that is, for every subset of nodes $S \subset V$ with $2 \leq |S| \leq |V| - 2$, the outflow is strictly greater than 1.

# 3 Vertex generation

In this section, we present our routine for the generation of half-integer vertices. First, making use of Theorem 2, we describe an efficient way of encoding cycle covers. We then discuss how the generated cycle covers are filtered to consider only non-isomorphic instances that satisfy the subtour elimination constraints (1d) and are extreme points of $\mathfrak{P}^n$. We say that two instances are isomorphic if the corresponding graphs are isomorphic as directed graphs. In what follows, we will also assume that $V^n = \{0, \dots, n-1\}$ for simplicity, yet similar results could be obtained for any given set of vertices $V^n$ by selecting an ordering of its elements.

## 3.1 Encoding of cycle covers

We first focus on the underlying cycle covers to establish a framework for their efficient generation. To begin, let us introduce a method for encoding each cycle cover, which allows us to reduce the problem of generating non-isomorphic cycle covers to the generation of their encodings (under certain clauses).

We encode any cycle cover by the sequences of vertices that each cycle in the cover visits, as in Figure 3. Each sequence starts from the lowest-numbered vertex and follows the arcs in the cycle to list the other vertices, while the sequences themselves are ordered by non-increasing cardinality (sequences with the same cardinality are ordered by increasing first element). This constructive arrangement guarantees that the encoding is unique for a given labeling of the vertices, removing redundant encodings due to cycle rotations. Definition 2 provides a formal definition of the cycle cover encoding.

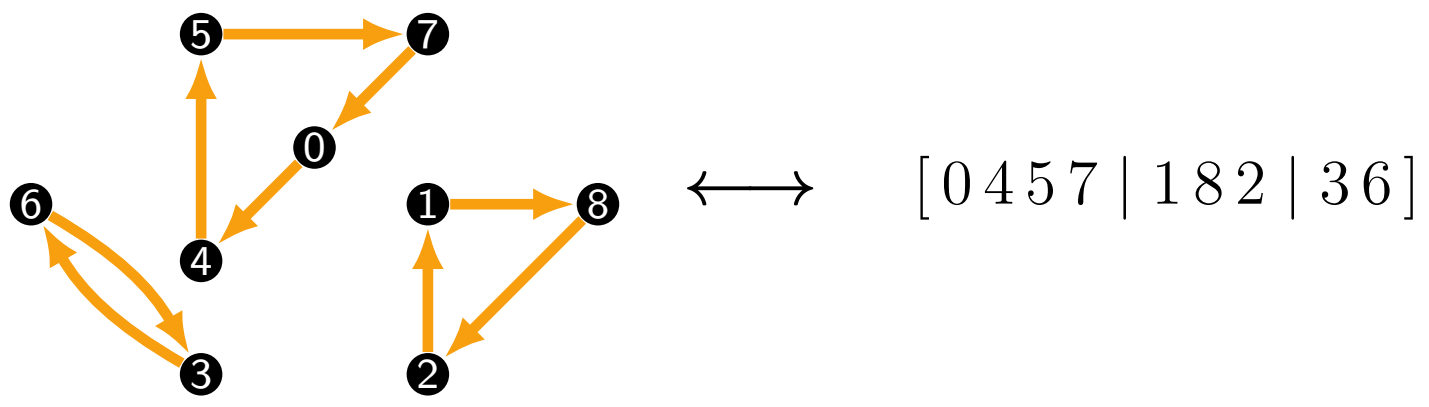

Figure 3: Example of a node cycle cover with the associated encoding. Referring to Definition 2, we have $n = 9$, $N = 3$, $p_1 = 4, p_2 = 3, p_3 = 2$.

**Definition 2** (Cover encoding)**.** Let $n, N \in \mathbb{N}$, $n \geq 2$, $1 \leq N \leq n$, and let $p_1, \dots, p_N$ be an integer partition of $n$, that is $p_1, \dots, p_N \in \mathbb{N}_{>0}$ such that $p_1 + \dots + p_N = n$ with

$p_1 \geq ... \geq p_N \geq 2$. We define

$$\xi = [\, k_1^1 \ldots k_{p_1}^1 \mid k_1^2 \ldots k_{p_2}^2 \mid \cdots \mid k_1^N \ldots k_{p_N}^N \,],$$

with $k_i^j \in V^n$ and $k_i^j \neq k_{i'}^{j'}$ for all $j, j', i, i' \in \mathbb{N}_{>0},\ j, j' \leq N,\ i \leq p_j,\ i' \leq p_{j'}$ as the encoding of the cover $C$ of $\mathcal{K}_n$ defined by

$$C \coloneqq \bigcup_{j=1}^{N} \left( \left\{ k_i^j k_{i+1}^j \,\middle|\, i \in \mathbb{N}_{>0}, i \leq p_j - 1 \right\} \cup \left\{ k_{p_j}^j k_1^j \right\} \right).$$

Let $(p_1, ..., p_N)$ be the *partition* of the encoding $\xi$. We also ask for the following two conditions to be met by the encoding:

i. $k_1^j < k_i^j$ for all $1 \leq j \leq N,\ 2 \leq i \leq p_j$;

ii. $k_1^j < k_1^{j'}$ for any $j < j'$ such that $p_j = p_{j'}$.

These last two requirements are necessary to pinpoint a specific encoding among all those that return the same cover $C$. Notice that any cover $C$ of a graph on $V^n$ has an associated encoding, which can be obtained by following Algorithm 1 reported in Appendix A: the selection of the cycles of $C$ in an ordered manner ensures that the encoding requirements (i.) and (ii.) are met. The existence of the algorithm also implies that for a fixed vertex set $V^n$, the encoding is unique.

Given a permutation $\sigma\colon \{0, \ldots, n-1\} \hookrightarrow \{0, \ldots, n-1\}$ and an encoding $\xi$, we write $\sigma(\xi)$ to identify the *translation* of the encoding via $\sigma$, that is the element-wise application of $\sigma$ to $\xi$ together with an eventual reordering of its elements so that the resulting encoding satisfies both requirements (i.) and (ii.) of Definition 2. The cover encoded by the translation $\sigma(\xi)$ is trivially isomorphic to that encoded by $\xi$, as the permutation $\sigma$ itself is an arc-preserving isomorphism between the two. The mentioned reordering is equivalent to and can be implemented as the execution of Algorithm 1 to the result of the application of $\sigma$ to $\xi$.

If we relax the requirements of Definition 2 by letting the partition of the encoding be the integer partition of an $m \leq n$, we are able to define the concatenation $\xi_1 + \xi_2$ of two such encoding-slices $\xi_1$ and $\xi_2$, with $\sum_{j=1}^{M'} p_j = m', \sum_{j=1}^{M''} q_j = m''$ as partitions respectively and $m' + m'' \leq n$, as

$$\xi_1 = [k_1^1 \ldots k_{p_1}^1 \mid \cdots \mid k_1^{M'} \ldots k_{p_{M'}}^{M'}], \qquad \xi_2 = [h_1^1 \ldots h_{q_1}^1 \mid \cdots \mid h_1^{M''} \ldots h_{q_{M''}}^{M''}],$$
$$\xi_1 + \xi_2 \coloneqq [k_1^1 \ldots k_{p_1}^1 \mid \cdots \mid k_1^{M'} \ldots k_{p_{M'}}^{M'} \mid h_1^1 \ldots h_{q_1}^1 \mid \cdots \mid h_1^{N''} \ldots h_{q_{M''}}^{M''}],$$

always ensuring that the resulting encoding-slice satisfies the other requirements of a cover encoding (e.g., $\xi_1$ and $\xi_2$ must have no common elements).

## 3.2 Encoding of sets of cycle covers

We now consider sets of cycle covers. We encode them by the encodings of their member covers as defined above, appropriately ordered with lexicographic criteria to make this new encoding unique (again, for a fixed labeling of the vertices). To remove part of the isomorphic instances that arise due to vertex relabeling, we further introduce the concept of *standard form*, obtained through a translation of the vertices that transforms the encoding of the first cover in the set to list all vertices in order. An example of cover-set encoding in standard form (that includes the cover set of Figure 3) can be found in Figure 4.

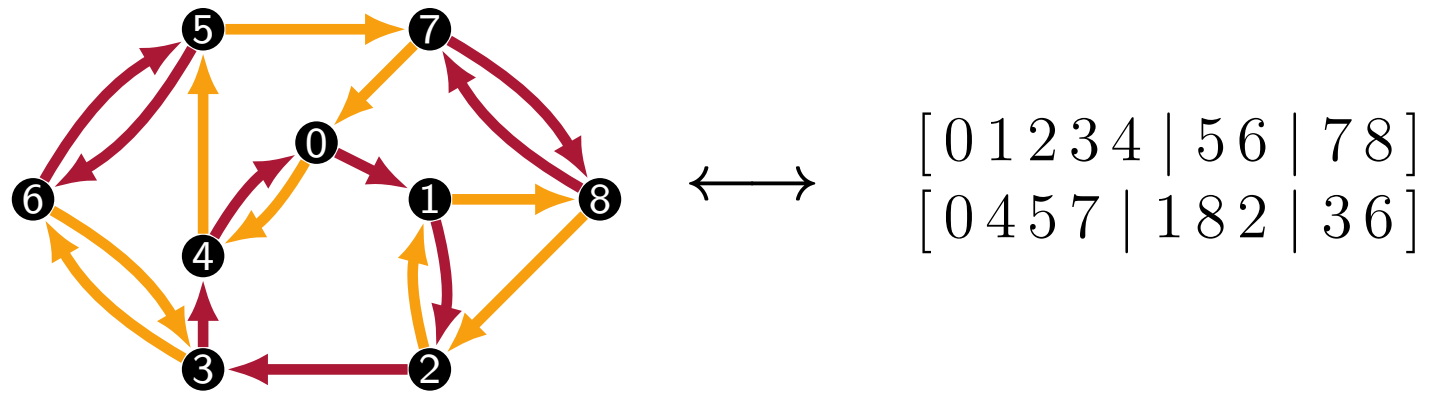

Figure 4: Example cover-set with its associated encoding (in standard form).

**Definition 3** (Cover-set encoding)**.** Let $\mathcal{C}$ be a set of cycle covers of $\mathcal{K}_n$. We encode $\mathcal{C}$ as the ordered set of encodings of its covers. The ordering is obtained by the following two criteria:

i. Descending lexicographic order of the encodings' partitions;

ii. For identical partitions, ascending lexicographic order of the encodings.

Note that when analyzing half-integer solutions, the cover sets always have cardinality 2, as they do throughout this work. More generally, one could consider points in $\mathfrak{X}^n_\alpha$ (see Remark 1), for which the corresponding cover sets would have cardinality $\alpha$.

While the cover-set encoding is unique when the vertex set $V^n$ is fixed (a direct consequence of the uniqueness of the cover encodings), this is not true for any relabeling of the vertices, a procedure that trivially preserves isomorphism. Therefore, it would be desirable to exploit isomorphism due to vertex relabeling to reduce the number of isomorphic instances considered by our generating procedure. We aim to find a permutation of the labels transforming the first of the cover encodings to the form

$$[\ 0 \ldots (p_1 - 1) \mid p_1 \ldots (p_1 + p_2 - 1) \mid \cdots \mid (n - p_N - 1) \ldots (n-1)\ ]. \tag{3}$$

We will refer to a cover-set encoding in this configuration as an encoding in *standard form*, see Algorithm 2 in Appendix A for a procedure to retrieve the standard form of any cover-set encoding. The existence of a standard form encoding isomorphic to any given encoding also proves that by generating all possible standard encodings, all cover-sets up to isomorphism are generated. In [24, S2.4], several algorithmic procedures are presented for generating all possible standard form encodings of cover-sets for a given $n$.

For a given set-cover, the standard form is however not unique, as in general the result of Algorithm 2 depends on the random selections made during its computation. Thus, we generate instances through their encoding in standard form, and then filter them for isomorphism.[1] Note that, when dealing with pure half-integer vertices only, we do not need to check for edge-weight preserving isomorphisms, since the weight of each arc is either 0 (no arc is present) or $\frac{1}{2}$. In our implementation, the isomorphism-class check relies on the C-based library Nauty [21]. See again [24, S2.3.2] for further details.

[1]The introduction of a *canonical form* of cover-set encodings (see Appendix B) would remove duplicate isomorphic encodings due to vertices relabeling. However, isomorphic set-covers with identical support yet different encodings in canonical form still exist. As the generation of canonical instances is furthermore computationally expensive, generating standard forms and later checking for isomorphism is therefore preferable.

## 3.3 Feasibility and extremality checking

Checking for the violation of the subtour elimination constraint is achieved by directly verifying the satisfaction by $\overline{\boldsymbol{x}}$ of inequalities (1d) by enumeration. We are then left with checking whether $\overline{\boldsymbol{x}}$ is also an extreme point of $\mathfrak{P}^n$ ($\overline{\boldsymbol{x}} \in \mathfrak{X}^n$). This is equivalent to asking for $\overline{\boldsymbol{x}}$ to be a basic solution, that is, a solution for which the number of linearly independent active constraints equals the dimension of the solution space.

Computing the matrix rank is a demanding task, and hence we elaborate a faster method to evaluate the rank-maximality of the active constraints. Such *circuit extremality algorithm* (Proposition 2) is extensively described in the following section.

# 4 Circuit Extremality Algorithm

In the present section, we develop a novel algorithm to fasten the check for the extremality of a given candidate half-integer solution. The designed algorithm relies on partitioning the non-zero arcs of $\overline{\boldsymbol{x}}$ based on the degree constraint vectors and its successive coarsening by the addition of the vectors associated with active subtour elimination constraints.

For pure half-integer vertices, the support of $\overline{\boldsymbol{x}}$ has two incoming and two outgoing arcs for each vertex in $V^n$, or equivalently that exactly two arcs are in the support of every degree constraint vector. The initial partitioning is built by repeatedly starting from an arc in the support not yet visited, and successively taking the one other arc having alternately the same head vertex or the same tail vertex as the previous one (i.e. the other arc in the support of the head's in-degree or tail's out-degree constraint vector), until the starting arc is reached. The arcs in each loop are alternately grouped in the same partition, thus forming two partitions (*circuits*) per loop that are *dual* to one another, and together form a *circuit pair*.

The circuit partitioning is then coarsened by successively taking each active subtour elimination constraint, and considering the two arcs in its support that are also in the support of $\overline{\boldsymbol{x}}$: if the two arcs belong to the same circuit, it is merged with its dual to a single partition, and we refer to the circuit and its pair as *shorted*; if the two arcs belong to two circuits in different pairs, then each of those circuits is merged with the dual of the other circuit to form a unique circuit pair (shorted if any of the two merged pairs is already shorted). This 'short-circuiting' procedure constitutes the core of the circuit extremality algorithm: the vertex $\overline{\boldsymbol{x}}$ is extremal if and only if at the end of the procedure the partitioning consists of a single shorted circuit.

**Theoretical Framework** We start by developing the necessary theory; the resulting *circuit extremality algorithm* is defined in Proposition 2. Proofs of the results stated in this section are provided in Appendix C.

If we analyze the coefficient vectors of the ASEP constraints, we see that they are all binary vectors in the space $\{0,1\}^m$, and they are fixed for any choice of $n$. Let us refer to $\deg_u^{+,n}, \deg_v^{-,n}$ as the vectors relative to the out-degree and in-degree constraints of $u, v \in V^n$ respectively, to $\mathrm{bnd}_{uv}^n$ as the vector of the bound constraint relative to $uv \in E^n$, and to $\mathrm{sep}_S^n$ as the vector of the subtour elimination constraint of $S \in \mathcal{S}$. Let us denote by $\mathcal{S}_{\mathtt{ACT}}^n(\overline{\boldsymbol{x}})$ the set of sets $S \in \mathcal{S}^n$ such that $\mathrm{sep}_S^n$ is the vector associated with an active constraint.

Thanks to the results proved in [15, S4.6], we can assume that $\overline{\boldsymbol{x}}$ is composed of only $0, 1/2$ entries ($\overline{\boldsymbol{x}} \in \mathfrak{P}_2^{n-1} \cap \{0, 1/2\}^m)^2$. Hence, we consider only the rank maximality over the subspace $\langle E^n_{\overline{\boldsymbol{x}}} \rangle := \operatorname{span}\left(\{\operatorname{bnd}^n_{uv} \mid uv \in E^n_{\overline{\boldsymbol{x}}}\}\right)$ with $E^n_{\overline{\boldsymbol{x}}} := \operatorname{supp}_{E^n}(\overline{\boldsymbol{x}})$. In the remainder, we will thus safely assume that all degree and active subtour elimination vectors are replaced with their projection $\lfloor \boldsymbol{a} \rfloor := \operatorname{pr}_{\langle E^n_{\overline{\boldsymbol{x}}} \rangle}(\boldsymbol{a})$ (by construction the active bound vectors are all projected to $\mathbf{0}$ and can be discarded).

Let us define the following sets of projected vectors:

$$\operatorname{Deg}^+_{\overline{\boldsymbol{x}}} := \{\lfloor \deg^{+,n}_u \rfloor \mid u \in V^n\}, \qquad \operatorname{Deg}_{\overline{\boldsymbol{x}}} := \operatorname{Deg}^+_{\overline{\boldsymbol{x}}} \cup \operatorname{Deg}^-_{\overline{\boldsymbol{x}}},$$
$$\operatorname{Deg}^-_{\overline{\boldsymbol{x}}} := \{\lfloor \deg^{-,n}_v \rfloor \mid v \in V^n\}, \qquad \operatorname{Sep}_{\overline{\boldsymbol{x}}} := \{\lfloor \operatorname{sep}^n_S \rfloor \mid S \in \mathcal{S}^n_{\mathtt{ACT}}(\overline{\boldsymbol{x}})\}.$$

Since $\forall \boldsymbol{a} \in \operatorname{Deg}_{\overline{\boldsymbol{x}}} \cup \operatorname{Sep}_{\overline{\boldsymbol{x}}}$ we have $\boldsymbol{a} \cdot \overline{\boldsymbol{x}} = 1$, considering that $\overline{x}_{uv} \in \{0, 1/2\}$ and $a_{uv} \in \{0, 1\}$ for all $uv \in E^n$, it follows that

$$|\operatorname{supp}_{E^n} \lfloor \boldsymbol{a} \rfloor| = |\operatorname{supp}_{E^n}(\boldsymbol{a}) \cap \operatorname{supp}_{E^n}(\overline{\boldsymbol{x}})| = |\operatorname{supp}_{E^n}(\boldsymbol{a} \cdot \overline{\boldsymbol{x}})| = \boldsymbol{a} \cdot 2\overline{\boldsymbol{x}} = 2,$$

showing how $\lfloor \boldsymbol{a} \rfloor$ is always a $0,1$-vector in $\langle E^n_{\overline{\boldsymbol{x}}} \rangle$ with exactly two 1-components. If $\operatorname{supp}_{E^n} \lfloor \boldsymbol{a} \rfloor = \{u_1v_1, u_2v_2\}$, we write $[u_1v_1|u_2v_2] = [u_2v_2|u_1v_1] := \lfloor \boldsymbol{a} \rfloor$ and we refer to $u_1v_1$ and $u_2v_2$ as the two *extremities* of $\lfloor \boldsymbol{a} \rfloor$ (notice also how all out- and in-degree vectors are respectively in the form $[uv_1|uv_2]$ and $[u_1v|u_2v]$). Finally, unraveling the definitions, we can write

$$[u_1v_1|u_2v_2] = \lfloor \operatorname{bnd}^n_{u_1v_1} \rfloor + \lfloor \operatorname{bnd}^n_{u_2v_2} \rfloor = \operatorname{bnd}^n_{u_1v_1} + \operatorname{bnd}^n_{u_2v_2}.$$

We now focus on the structure of $\operatorname{Deg}_{\overline{\boldsymbol{x}}}$: for any $uv \in E^n_{\overline{\boldsymbol{x}}}$ there are exactly two degree vectors with non-null $\overline{u}\,\overline{v}$-component, that are $\lfloor \deg^{+,n}_{\overline{u}} \rfloor = [\overline{u}\,\overline{v}|\overline{u}v]$ and $\lfloor \deg^{-,n}_{\overline{v}} \rfloor = [\overline{u}\,\overline{v}|u\overline{v}]$. This justifies the following definition.

**Definition 4** (Link, circuit, circuit partition)**.** For a set of constraints $A$ such that $\operatorname{Deg}_{\overline{\boldsymbol{x}}} \subseteq A \subseteq \operatorname{Deg}_{\overline{\boldsymbol{x}}} \cup \operatorname{Sep}_{\overline{\boldsymbol{x}}}$, let $\sim_A$ be the transitive closure of the homogeneous binary relation $\leftrightarrow_A$ on $E^n_{\overline{\boldsymbol{x}}}$ defined by $u_1v_1 \leftrightarrow_A u_2v_2$ if and only if

$$\exists \overline{u}\,\overline{v} \in E^n_{\overline{\boldsymbol{x}}} \quad \text{such that } [u_1v_1|\overline{u}\,\overline{v}], [\overline{u}\,\overline{v}|u_2v_2] \in A.$$

Additionally, $\leftrightarrow_A$ and $\sim_A$ are trivially reflexive and symmetrical, thus $\sim_A$ is an equivalence relation. For any $u_1v_1 \neq u_2v_2$, we refer to $\overline{u}\,\overline{v}$ as the *link* between $u_1v_1$ and $u_2v_2$ and we write $u_1v_1 \overset{\overline{u}\,\overline{v}}{\longleftrightarrow}_A u_2v_2$, and we refer to the set of arcs that are linked to others by $\overline{u}\,\overline{v}$ as $\operatorname{adj}_A(\overline{u}\,\overline{v}) := \{uv \in A \mid \exists u'v' \in A, uv \overset{\overline{u}\,\overline{v}}{\longleftrightarrow}_A u'v'\}$.

Let the *circuit partition* $\mathcal{L}_A$ be the partitioning of the arcs in $E^n_{\overline{\boldsymbol{x}}}$ given by the equivalence classes of $\sim_A$, and call *circuits* its elements $L \in \mathcal{L}_A$. Since $\forall u_1v_1, u_2v_2, u_3v_3 \in L$ we have

$$u_1v_1 \overset{\overline{u}_1\overline{v}_1}{\longleftrightarrow}_A u_2v_2 \overset{\overline{u}_2\overline{v}_2}{\longleftrightarrow}_A u_3v_3 \quad \Longrightarrow \quad \overline{u}_1\overline{v}_1 \overset{u_2v_2}{\longleftrightarrow}_A \overline{u}_2\overline{v}_2,$$

it follows that there exists a $L' \in \mathcal{L}_A$ containing both $\overline{u}_1\overline{v}_1$ and $\overline{u}_2\overline{v}_2$. We refer to such $L'$ as the *dual* of $L$ denoted as $\overline{L}^A$, and refer to the set $\{L, \overline{L}^A\}$ as a *circuit pair*. We say that a $L \in \mathcal{L}_A$ for which $\overline{L}^A = L$ is *short-circuited* or simply *shorted*.

Note that $A \subseteq A'$ implies $\leftrightarrow_A \subseteq \leftrightarrow_{A'}$ and $\sim_A \subseteq \sim_{A'}$ and thus makes $\mathcal{L}_{A'}$ a coarser partition than $\mathcal{L}_A$. All the subscripts $A$ of $\leftrightarrow_A$, $\sim_A$, $\operatorname{adj}_A$, $\mathcal{L}_A$ $\overline{\cdot}^A$ will be dropped whenever it is unambiguous to do so.

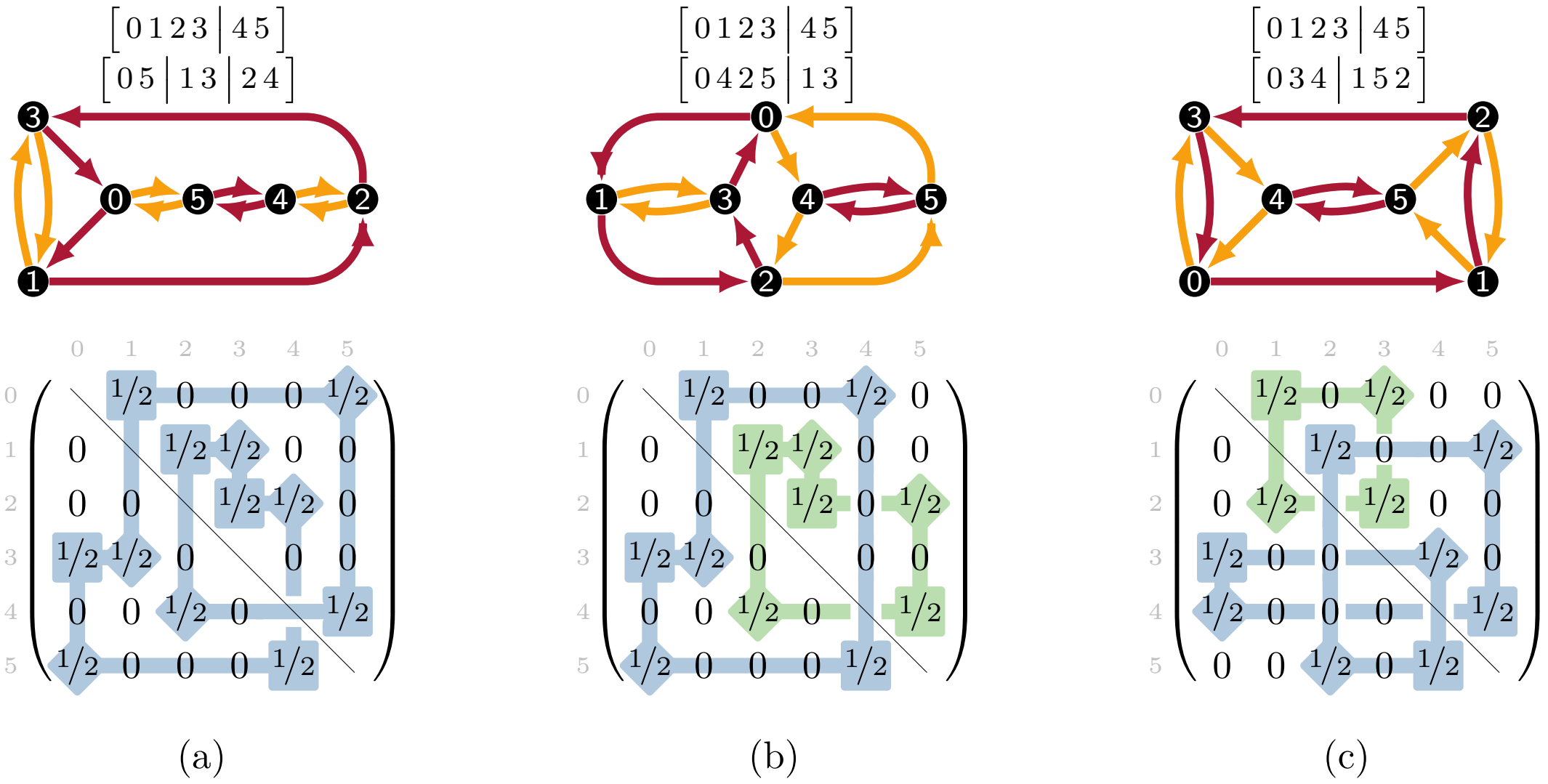

Figure 5: Examples of circuit partitions for $A = \mathrm{Deg}_{\overline{\boldsymbol{x}}}$ on the matrix representations of the solution's characteristic vectors.

**Remark 2.** The most interesting case is $A = \mathrm{Deg}_{\overline{\boldsymbol{x}}}$. Again, since any $\overline{u}\,\overline{v} \in E^{n}_{\overline{\boldsymbol{x}}}$ is the extremity of an out-degree $[\overline{u}\,\overline{v}|\overline{u}v]$ and an in-degree $[\overline{u}\,\overline{v}|u\overline{v}]$ vector we have that every $\overline{u}\,\overline{v} \in E^{n}_{\overline{\boldsymbol{x}}}$:

- ▷ it links two distinct arcs $\overline{u}v, u\overline{v}$ via an out-degree and an in-degree vector;
- ▷ it is linked to two (eventually equal) arcs $u'v, uv'$ by $\overline{u}v, u\overline{v}$, respectively.

Hence, any circuit pair $\left\{L, \overline{L}\right\}$ identifies a succession of linked arcs, which eventually loops on itself (due to the finiteness of the set of arcs)

$$\circlearrowright\; u_1v_1 \xleftrightarrow{\overline{u}_1\overline{v}_1} u_2v_2 \xleftrightarrow{\overline{u}_2\overline{v}_2} u_3v_3 \longleftarrow \cdots \longrightarrow u_{|L|}v_{|L|} \xleftarrow{\overline{u}_{|L|}\overline{v}_{|L|}}$$

with $u_iv_i \in L$, $\overline{u}_i\overline{v}_i \in \overline{L}$ for all $i = \{1, ..., |L|\}$. This property shows that $|L| = \left|\overline{L}\right|$. By taking into account the resulting "domino" succession of vectors

$$\Big( [u_1v_1|\overline{u}_1\overline{v}_1]\,, [\overline{u}_1\overline{v}_1|u_2v_2]\,, \ldots\,, [u_{|L|}v_{|L|}|\overline{u}_{|L|}\overline{v}_{|L|}]\,, [\overline{u}_{|L|}\overline{v}_{|L|}|u_1v_1] \Big), \tag{4}$$

it follows that it must be composed of alternating out- and in-degree constraint vectors. But this, in turn implies the non-trivial definition of the dual of a circuit: we have that $L$ is not short-circuited ($L \neq \overline{L}$), since if it were $u_iv_i = \overline{u}_j\overline{v}_j$ for some $i, j$ we would have either $[u_iv_i|\overline{u}_i\overline{v}_i]\,, [u_jv_j|\overline{u}_j\overline{v}_j]$ or $[\overline{u}_j\overline{v}_j|u_{j+1}v_{j+1}]\,, [\overline{u}_{i-1}\overline{v}_{i-1}|u_iv_i]$ (where the subscript is adjusted modulo $|L|$) as a pair of both out- or in-degree vectors with the same arc as extremity, against the hypothesis. See Figure 5 for an example. There, circuit pairs are displayed as closed curves of different colors, with duals represented as square/diamond nodes. Two entries of the matrix are linked horizontally if the two corresponding edges guarantee that the out-degree constraint of the node corresponding to the row is satisfied, while two entries of the matrix are linked vertically if the same happens for the in-degree constraint of the node corresponding to the column. For example, in the first matrix, the entries $(0,1)$ and $(0,5)$ are linked horizontally because they are the outgoing arcs of the node 0, while the entries $(0,1)$ and $(3,1)$ are linked vertically because they are the incoming arcs of the node 1.

[2] The circuit extremality algorithm we introduce is however also valid on general half-integer vertices, as explained in [24, S2.3.3].

For $L, \overline{L} \in \mathcal{L}_A$, let us denote by $[L]_A$ the subset of constraints of $A$ that constitute the linking of $L$ and $\overline{L}$, formally $[L]_A \coloneqq \left\{ [uv|\overline{u}\,\overline{v}] \in A \;\middle|\; uv \in L, \overline{u}\,\overline{v} \in \overline{L} \right\}$. With regards to the subspaces $\operatorname{span}[L]_A$ generated by the circuits, we notice that distinct circuit pairs generate linearly independent subspaces, as the vectors $[L]_A$ have extremities (and thus support) entirely contained in $L \cup \overline{L}$ and the circuits themselves are a partition of $E^n_{\overline{\boldsymbol{x}}}$. Considering again the case of $A = \operatorname{Deg}_{\overline{\boldsymbol{x}}}$, we can see that $[L]_{\operatorname{Deg}_{\overline{\boldsymbol{x}}}}$ is equivalent to the set expressed in (4), and we can thus decompose it as

$$[L]_{\operatorname{Deg}_{\overline{\boldsymbol{x}}}} = [L]^+ \cup [L]^- , \qquad \text{with} \begin{cases} [L]^+ \coloneqq [L]_{\operatorname{Deg}_{\overline{\boldsymbol{x}}}} \cap \operatorname{Deg}^+_{\overline{\boldsymbol{x}}}, \\ [L]^- \coloneqq [L]_{\operatorname{Deg}_{\overline{\boldsymbol{x}}}} \cap \operatorname{Deg}^-_{\overline{\boldsymbol{x}}}. \end{cases}$$

Under this assumption, we can observe that $\operatorname{span}[L]_{\operatorname{Deg}_{\overline{\boldsymbol{x}}}}$ has dimension $2|L| - 1$, since $\dim(\operatorname{span}[L]^+) = \dim(\operatorname{span}[L]^-) = |L|$ (both $[L]^+$ and $[L]^-$ are sets of linearly independent vectors) and $\dim(\operatorname{span}([L]^+) \cap \operatorname{span}([L]^-)) = 1$ from

$$\begin{aligned} \operatorname{span}\left([L]^+\right) \cap \operatorname{span}\left([L]^-\right) &= \operatorname{span}\left\{ \textstyle\sum_{[e_1|e_2] \in [L]^+} [e_1|e_2] \right\} = \\ &= \operatorname{span}\left\{ \textstyle\sum_{[e_2'|e_1'] \in [L]^-} [e_2'|e_1'] \right\} = \operatorname{span}\left\{ \textstyle\sum_{e \in L \cup \overline{L}} \operatorname{bnd}^n_e \right\}. \end{aligned}$$

Therefore, we aim at finding a subtour elimination constraint in $\operatorname{Sep}_{\overline{\boldsymbol{x}}}$ that "fills" the remaining degree of freedom of $\operatorname{span}[L]_{\operatorname{Deg}_{\overline{\boldsymbol{x}}}}$. A characterization of those "filling" constraints is given by the following propositions.

**Proposition 1** (Circuits merging)**.** *Let* $[e_1|e_2] \in \operatorname{Sep}_{\overline{\boldsymbol{x}}}$ *and define* $A' \coloneqq A \cup \{[e_1|e_2]\}$ *for a given* $A$ *as in Definition 4. Given* $L_1, L_2 \in \mathcal{L}_A$ *such that* $e_1 \in L_1, e_2 \in L_2$*, then*

$$\mathcal{L}_{A'} = \mathcal{L}' \cup \mathcal{L}_A \setminus \{L_1, L_2, \overline{L_1}^A, \overline{L_2}^A\},$$

*with*

*i.* $\mathcal{L}' \coloneqq \{L_1 \cup \overline{L_1}^A \cup L_2 \cup \overline{L_2}^A\}$ *if* $L_1 = L_2$*,* $L_1 = \overline{L_1}^A$ *or* $L_2 = \overline{L_2}^A$*, where the only circuit in* $\mathcal{L}'$ *is shorted.*

*ii.* $\mathcal{L}' = \{L_1 \cup \overline{L_2}^A, \overline{L_1}^A \cup L_2\}$ *otherwise, where the two circuits in* $\mathcal{L}'$ *are duals and not shorted.*

*In addition, if* $L_2 = \overline{L_1}^A$*, then* $\mathcal{L}_{A'} = \mathcal{L}_A$*.*

**Corollary 1.** *In the above notation, if* $L_2 = \overline{L_1}^A$*, then* $\mathcal{L}_{A'} = \mathcal{L}_A$*.*

**Lemma 2.** *Let* $L \in \mathcal{L}_A$ *such that* $\dim(\operatorname{span}[L]_A) = 2|L| - 1$ *and* $[e_1|e_2] \in \operatorname{Sep}_{\overline{\boldsymbol{x}}} \setminus A$ *such that* $e_1 \in L$*. Then* $[e_1|e_2]$ *is linearly independent from* $[L]_A$ *if and only if* $e_2 \notin \overline{L}^A$*.*

**Lemma 3.** *Let* $L \in \mathcal{L}_A$*. Then* $\dim(\operatorname{span}[L]_A)$ *is equal to* $|L|$ *if* $L$ *is shorted and to* $2|L| - 1$ *otherwise.*

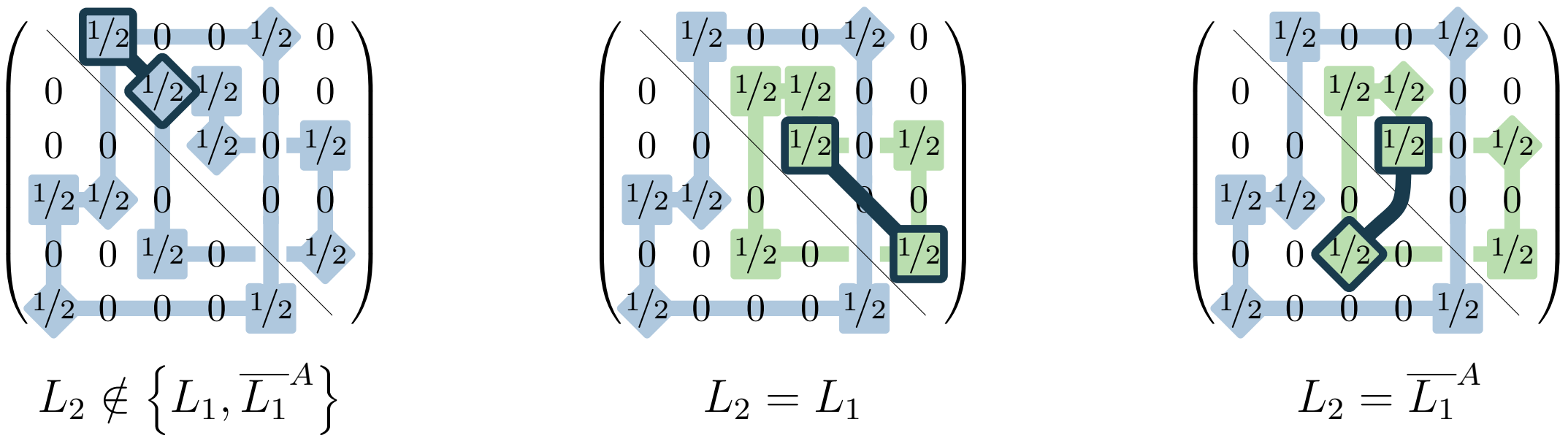

Figure 6: Schematic representation of different possible circuits merging in the example of Figure 5b; in the notation of Proposition 1, the dark-colored connections represents the newly added vector $[e_1|e_2]$ with $e_1 \in L_1, e_2 \in L_2$, while the partitions resulting from the merging are represented similarly to Figure 5, with shorted partitions shown having only square or only diamond nodes.

**Circuit Extremality Algorithm** Having provided all the required results, we can now define the procedure to check for the extremality of a given solution $\overline{\boldsymbol{x}} \in \mathfrak{P}_2^n$, thus proving whether $\overline{\boldsymbol{x}} \in \mathfrak{X}_2^n$. The algorithm considers the arcs $e \in E^n$ for which $\overline{x}_e = 1/2$ and partitions them in circuits using the set of degree vectors; the partitions are then repeatedly merged using the vectors of active subtour elimination constraints. The algorithm stops whenever the circuit pairs in the partitioning are all shorted, which guarantees the rank maximality of the constraints active at $\overline{\boldsymbol{x}}$ and thus, the extremality of the solution. If the subtour elimination constraints are exhausted before the terminating condition is met, the active constraints are not maximal and $\overline{\boldsymbol{x}} \notin \mathfrak{X}_2^n$. Checking for the terminating condition is not achieved in a direct manner, but by counting the times when a partition is shortened by merging it with its dual or by merging it with another shortened partition. The full procedure is detailed in Algorithm 3 of Appendix A. Compared to our first nontrivial implementation, the circuit extremality algorithm significantly improves the runtime, with a 350× speedup for $n = 11$.

**Proposition 2.** *Given $\overline{\boldsymbol{x}} \in \mathfrak{P}_2^n$, proceed as follows.*

*i. Find $E_{\overline{\boldsymbol{x}}}^n = \{e \in E^n \mid \overline{x}_e = 1/2\}$ and determine*

$$\mathrm{Deg}_{\overline{\boldsymbol{x}}} \coloneqq \left\{ \mathrm{pr}_{\langle E_{\overline{\boldsymbol{x}}}^n \rangle} \left(\mathrm{deg}_u^{\pm,n}\right) \;\middle|\; u \in V^n \right\}, \quad \mathrm{Sep}_{\overline{\boldsymbol{x}}} \coloneqq \left\{ \mathrm{pr}_{\langle E_{\overline{\boldsymbol{x}}}^n \rangle} \left(\mathrm{sep}_S^n\right) \;\middle|\; S \in \mathcal{S}_{\mathtt{ACT}}^n(\overline{\boldsymbol{x}}) \right\},$$

*and compute $\mathcal{L}_A$ having initialized $A = \mathrm{Deg}_{\overline{\boldsymbol{x}}}$, then set $c = \frac{1}{2}|\mathcal{L}_A|$.*

*ii. Select $[e_1|e_2] \in \mathrm{Sep}_{\overline{\boldsymbol{x}}} \setminus A$ and add it to $A$, merging the partitions of $L_1, L_2 \in \mathcal{L}_A$ as in Proposition 1: if $L_2 \neq \overline{L_1}$ and not both $L_1, L_2$ are shorted, decrease $c$ by 1.*

*iii. If $c$ has reached 0, then $\overline{\boldsymbol{x}} \in \mathfrak{X}_2^n$, otherwise repeat point* (ii.)*. If all elements of $\mathrm{Sep}_{\overline{\boldsymbol{x}}}$ have been selected, then $\overline{\boldsymbol{x}} \notin \mathfrak{X}_2^n$.*

When implementing this algorithm, it turns out to be more efficient to assign to each arc a label corresponding to the circuit the arc belongs to. This allows for a faster retrieval of the partitions $L_1, L_2$ that contain the extremities of the selected active subtour elimination vector at each iteration. However, updating the arc labels at each partition merging becomes highly impractical. Therefore, it is preferable to leave unchanged the label of the initial partitioning at each arc, separately keeping track of which partitions get merged together and of the duals of each partition.

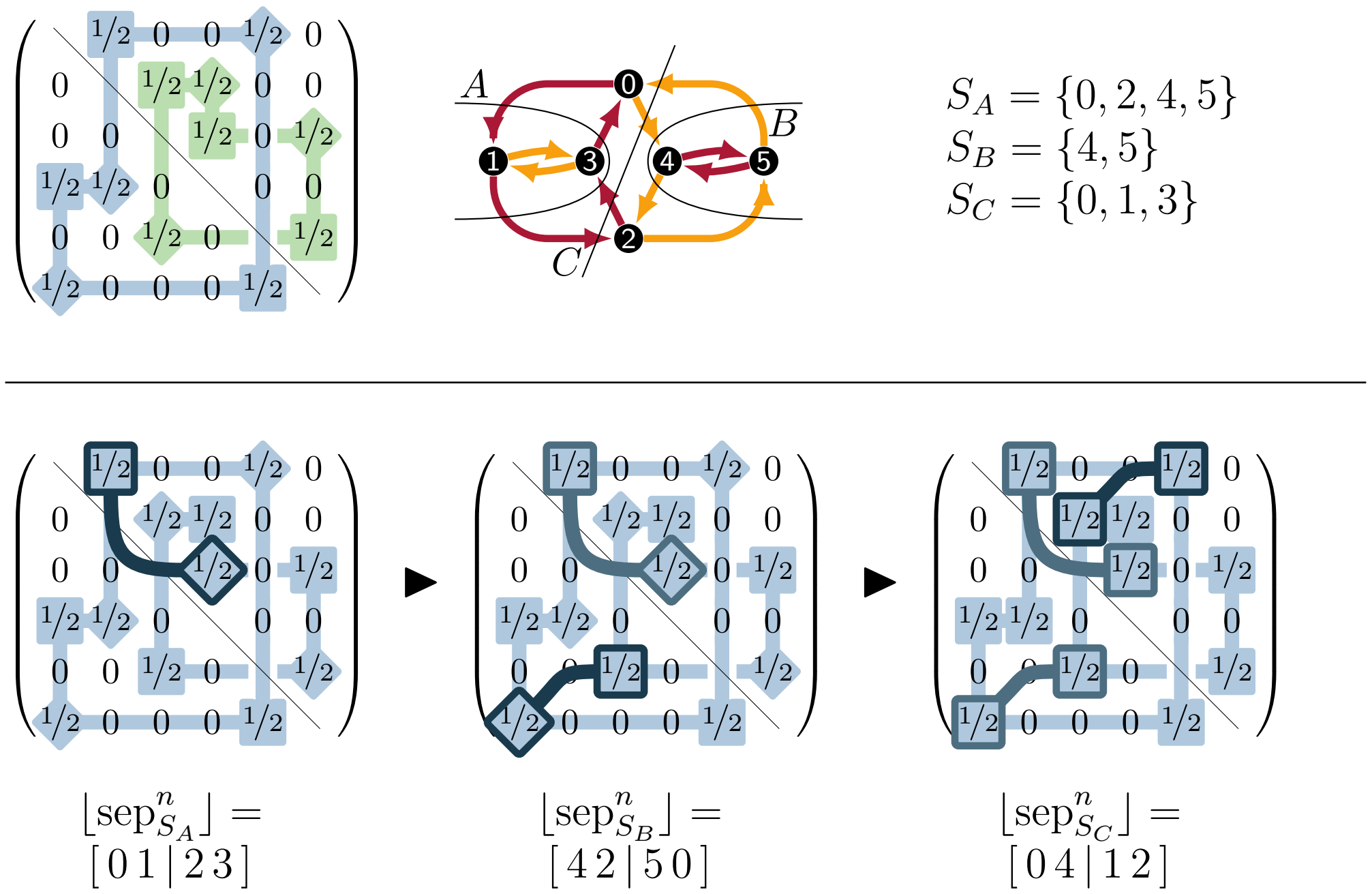

Figure 7: Schematic representation of a circuit extremality algorithm example run over the instance $\overline{\boldsymbol{x}}$ of Figure 5b; the projections of the active subtour elimination vectors relative to the subsets $S_A, S_B, S_C$ are successively added, and at each step the effects on the circuit partitions are shown similarly to the previous figures: the existence of a single, shorted circuit at the end proves that the input $\overline{\boldsymbol{x}}$ is an extreme point of $\mathfrak{P}^n$.

# 5 Computational results

We have implemented the algorithm described in Figure 1 mainly in Python, obtaining the results reported in Table 1, while vertices with the highest integrality gap are presented in Appendix D. Our results considerably improve over the ones presented in [15]: we replicated their results in a tiny fraction of time for any $n \leq 9$, and we extended the exact computation of half-integers for $n = 10, 11, 12$ in around 12 hours, proving computationally that the best-known bounds for the integrality gap up to $n = 12$ are indeed optimal. For $n \geq 9$, the maximum gap is attained by instances belonging to the recursive family introduced by Elliott-Magwood [15, S5.1], obtained via repeated applications of the 1-arc insertion and 2-jack operations starting from the $n = 4$ extremal example, and by a unique instance for $n = 9, 11, 12$, as depicted in Appendix D. Unfortunately, the author reports no runtimes except that the exhaustive ASEP vertices enumeration of $n = 7$ took over 20 hours of computation.

It is important to note that Table 1 highlights the rapid combinatorial growth of the enumeration process. The number of generated instances increases by roughly one order of magnitude when passing from $n$ to $n + 1$, reflecting the explosion in the number of cycle covers and their pairwise combinations. Since each candidate half-integer vertex arises from a pair of cycle covers, the size of the search space is governed by the combinatorics of cycle-cover structures, whose number grows rapidly with $n$. Although only a subset of such pairs yields actual extreme points, and isomorphism filtering drastically reduces redundancy by roughly two orders of magnitude, the empirical data indicate that this structural filtering does not significantly alter the overall growth trend. In particular, the number of non-isomorphic half-integer vertices grows similarly to that of the underlying cycle-cover constructions, suggesting that the half-integer vertices are closely tied to the combinatorial structure of cycle covers. From a computational standpoint, all phases of the pipeline exhibit exponential growth. For $n = 12$, the subtour and extremality verification stage dominates the overall runtime. However, the generation phase represents an inherent scalability barrier: even before extremality and gap computations are performed, the sheer number of generated candidates becomes prohibitive for larger $n$.

Table 1: Exact integrality gaps of pure half-integer vertices, number of generated and non-isomorphic vertices, and runtime for the subprocedures of our enumeration algorithm.

| $n$ | $\text{Gap}_n$ | Number of instances generated | Number of non-isom. vertices | Runtime GENER.+ CERTIFIC. | Runtime SUBT.+ EXTREM. | Runtime INTEGR. GAP | Runtime TOTAL |
|---|---|---|---|---|---|---|---|
| 5 | 5/4 | 11 | 2 | 0.22s | 0.26s | 0.28s | 0.76s |
| 6 | 4/3 | 112 | 11 | 0.30s | 0.28s | 1.30s | 1.90s |
| 7 | 4/3 | 978 | 52 | 0.21s | 0.32s | 1.30s | 1.80s |
| 8 | 4/3 | 10 809 | 365 | 0.35s | 0.55s | 1.50s | 2.40s |
| 9 | 11/8 | 116 463 | 2 931 | 1.60s | 0.78s | 7.90s | 10.00s |
| 10 | **7/5** | 1 527 602 | 26 906 | 17.00s | 1.80s | 53.00s | 1m 11s |
| 11 | **10/7** | 19 938 367 | 274 134 | 4m 34s | 13.00s | 14m 01s | 18m 48s |
| 12 | **56/39** | 306 394 067 | 3 059 487 | 1h 34m | 7h 12m | 4h 05m | 12h 51m |

**Implementation details.** The isomorphism-class check relies on the C-based library Nauty [21] via the Python package pynauty [14], while the remaining property-checking routines are implemented directly in Python. The state-of-the-art commercial solver Gurobi [18] is used for solving the GAP problem for each non-isomorphic vertex via the `gurobipy` interface. The software ran on a Dell Precision 7960 Tower Workstation, which features an Intel Xeon W-3495X CPU with 56 physical cores and 128 GB of memory. Chapter 3 of [24] provides a more detailed description of our experimental setup, further developed for the present work, together with an in-depth analysis of runtimes and resources required by each structural component of our pipeline and a discussion of the different strategies and configurations investigated.

**Data Availability Statement.** The code developed for this study, which supports the results and analysis presented in the article, is available at [25]. No additional datasets were generated or analyzed.

# 6 Conclusions and future work

In this paper, we have presented a novel approach to compute the integrality gap of small ATSP instances that enabled the computation of the exact integrality gap for half-integer vertices for $n = 10, 11, 12$. In our algorithm, we first generate pairs of edge-disjoint cycle covers through a lexicographically ordered encoding that filters out several isomorphic instances, and later, we efficiently check a few properties that characterize those cycle covers as a vertex of ASEP. In future works, we plan to apply our approach to larger values of $n$ and to $k$-tuple of cycle covers, $k \geq 3$, to study different classes of fractional vertices, and to adapt the described approach to the symmetric TSP.

# Appendix A Algorithms

In this appendix, we report the pseudocode of the algorithm described in the paper. The corresponding implementation in Python is available at [25].

**Algorithm 1** Cover encoding

▷ **Input** A cycle cover $C$
▷ **Output** A cover encoding $[k_1^1 \dots k_{p_1}^1 \mid k_1^2 \dots k_{p_2}^2 \mid \cdots \mid k_1^N \dots k_{p_N}^N]$

Sort the cycles in $C$ by decreasing cardinality
Sort the same-cardinality cycles in $C$ by increasing minimal vertex
$j \leftarrow 1$
**for** $c \in$ (sorted set of cycles of C) **do**
    $p_j \leftarrow |c|$
    $k_1^j \leftarrow$ minimal vertex in $c$
    **for** $i \in \{2, \dots, p_j\}$ **do**
        Find $uv \in \delta^+(k_{i-1}^j) \cap C$
        $k_i^j \leftarrow v$
    **end for**
    $j \leftarrow j + 1$
**end for**

**Algorithm 2** Standard cover-set encoding

▷ **Input** A cover-set encoding $\Xi = (\xi_1, \ldots, \xi_M)$, $(p_1, ..., p_N)$ partition of $\xi_1$

▷ **Output** A cover-set encoding $\Xi'$ in standard form

Let $\sigma\colon \{0, \ldots, n-1\} \hookrightarrow \{0, \ldots, n-1\}$ be the translating permutation
Select $\xi_m \in \Xi$ with the same partition $(p_1, ..., p_N)$ as $\xi_1$
Let $\xi_m = [k^1_1 \ldots k^1_{p_1} \mid \cdots \mid k^N_1 \ldots k^N_{p_N}]$
$P \leftarrow \emptyset$; $S \leftarrow 0$
**for** $j' \in \{1, \ldots, N\}$ **do**
    Select $j$ such that $p_j = p_{j'}$ and $j \notin P$
    $P \leftarrow P \cup \{j\}$
    Select $s \in \{0, \ldots, p_j - 1\}$
    **for** $i \in \{1, \ldots, p_j\}$ **do**
        $t \leftarrow s + i \mod p_j$
        $\sigma(k^j_i) \leftarrow S + t$
    **end for**
    $S \leftarrow S + p_j$
**end for**
$\Xi' \leftarrow (\sigma(\xi_1), \ldots, \sigma(\xi_M))$
Reorder $\Xi'$ following criteria (i.) and (ii.) of Definition 3

# Appendix B Canonical Form of a Cover-set Encoding

The standard form of a cover-set encoding is not unique. One may define the *canonical form* of a cover-set encoding the standard form that achieves the minimum in a lexicographic order among all possible standard forms (the comparison between elements of the standard forms is again done lexicographically). By construction, this canonical form is therefore unique for each set-cover encoding. This construction presents two issues. First, producing canonical instances indeed requires exploring all their possible standard translations to check whether the minimum has been reached, obviously an extremely demanding operation. Notice also how it is, in fact, equivalent to generating all standard forms and then checking for their canonicity afterward. Second, counterexamples, as shown in Figure 8, display two different cover sets (even with different partitions) in canonical form that represent two isomorphic vertices.

The canonical form is therefore not strong enough to identify the isomorphism class of a cover-set. It is therefore equivalent and better for us to generate all standard encodings, since it leaves us the choice of whether to check afterward for their canonicity or directly for their isomorphism class with faster algorithms.

In particular, it is computationally more efficient to just check if a generated instance is isomorphic to an instance we already encountered.

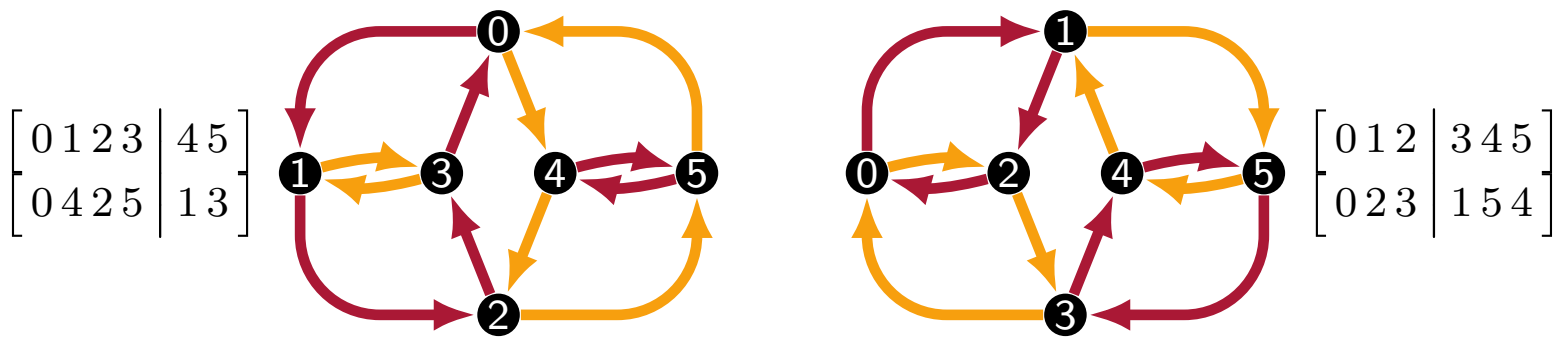

Figure 8: Isomorphic vertices with different cover-set canonical encodings. Note that the isomorphism sending the first one into the second one is represented by $\sigma := (015)(234)$.

**Algorithm 3** Circuit extremality algorithm

▷ **Input** A solution $\overline{\boldsymbol{x}} \in \mathfrak{P}_2^n$

▷ **Output** The membership of $\overline{\boldsymbol{x}}$ in $\mathfrak{X}_2^n$

**function** MergePartitions($L_1, L_2, \mathsf{duals}, \mathsf{shorted}$)
  **if** $L_1 \neq L_2 \ \wedge\ \{L_1, L_2\} \cap \mathsf{shorted} \neq \emptyset$ **then**
    Get $L \in \{L_1, L_2\}$ with $L \notin \mathsf{shorted}$
    $\mathsf{duals}\,(L) \leftarrow \mathsf{duals}\,(L) \cup \{L\}$
  **end if**
  **for** $(A, B) \in \{(L_1, L_2), (L_2, L_1)\}$ **do**
    $\mathsf{duals}\,(A) \leftarrow \mathsf{duals}\,(A) \cup \bigcup_{L \in \mathsf{duals}(B)} \mathsf{duals}\,(L)$
  **end for**
  **for** $(A, B) \in \{(L_1, L_2), (L_2, L_1)\}$ **do**
    **for** $L \in \mathsf{shorted}\,(B)$ **do**
      $\mathsf{duals}\,(L) \leftarrow \mathsf{duals}\,(A)$
    **end for**
    **if** $A \in \mathsf{duals}\,(A)$ **then**
      $\mathsf{shorted} \leftarrow \mathsf{shorted} \cup \mathsf{duals}\,(A)$
    **end if**
  **end for**
  **return** $\mathsf{duals}$, $\mathsf{shorted}$
**end function**

**function** CircuitPartition($\overline{\boldsymbol{x}}$)
  $E^n_{\overline{\boldsymbol{x}}} \leftarrow \{e \in E^n \mid \overline{x}_e = 1/2\}$
  $c \leftarrow 0$; $E \leftarrow E^n_{\overline{\boldsymbol{x}}}$; $\mathsf{part}\colon E^n_{\overline{\boldsymbol{x}}} \to \mathbb{N} \times \{0, 1\}$
  **repeat**
    $c \leftarrow c + 1$
    Select $uv \in E$
    **repeat**
      $\mathsf{part}\,(uv) \leftarrow (c, 0)$; $E \leftarrow E \setminus \{uv\}$
      Select $u'v'$ from $E^n_{\overline{\boldsymbol{x}}}$ with $u' = u, v' \neq v$
      $\mathsf{part}\,(uv) \leftarrow (c, 1)$; $E \leftarrow E \setminus \{u'v'\}$
      Select $uv$ from $E^n_{\overline{\boldsymbol{x}}}$ with $u \neq u', v = v'$
    **until** $uv \in E$
  **until** $E = \emptyset$
  **return** $c$, $\mathsf{part}$
**end function**

$c$, $\mathsf{part} \leftarrow$ CircuitPartition $(\overline{\boldsymbol{x}})$
$\mathsf{duals}\colon (k, i) \mapsto \{(k, 1 - i)\}$; $\mathsf{shorted} \leftarrow \emptyset$
**for** $[e_1, e_2] \in \mathrm{Sep}_{\overline{\boldsymbol{x}}}$ **do**
  $L_1 \leftarrow \mathsf{part}\,(e_1)$; $L_2 \leftarrow \mathsf{part}\,(e_2)$
  **if** $L_2 \notin \mathsf{duals}\,(L_1) \ \wedge\ (L_1 \notin \mathsf{shorted} \ \vee\ L_2 \notin \mathsf{shorted})$ **then**
    $c \leftarrow c - 1$
    **if** $c = 0$ **then**
      **return** `True`
    **end if**
    $\mathsf{duals}$, $\mathsf{shorted} \leftarrow$ MergePartitions $(L_1, L_2, \mathsf{duals}, \mathsf{shorted})$
  **end if**
**end for**
**return** `False`

# Appendix C Omitted Proofs

*Proof of Proposition 1.* We already remarked how $\leftrightarrow_A \subseteq \leftrightarrow_{A'}$, and we now construct explicitly the elements in $\leftrightarrow_{A'} \setminus \leftrightarrow_A$. We observe that the element $[e_1|e_2]$ can only be involved in the linkage of one of its extremities to arcs adjacent (with respect to $A$) to the other one, that is

$$e_1 \overset{[e_1|e_2]}{\longleftrightarrow}_{A'} \overline{e}_2 \quad \forall \overline{e}_2 \in \mathrm{adj}_A\,(e_2) \qquad \text{and} \qquad e_2 \overset{[e_1|e_2]}{\longleftrightarrow}_{A'} \overline{e}_1 \quad \forall \overline{e}_1 \in \mathrm{adj}_A\,(e_1)\,.$$

Since $\mathrm{adj}_A(e_2) \neq \emptyset$ (Remark 2) and $\mathrm{adj}_A(e_2) \subseteq \overline{L_2}^A$, the two equivalence classes $L_1$ and $\overline{L_2}^A$ having some $A'$-linked representatives thus merge under the new equivalence relation $\sim_{A'}$. Similarly happens for $L_2$ and $\overline{L_1}^A$. This means that the four (eventually non-distinct) $\sim_A$-equivalence classes merge into the $\sim_{A'}$-equivalence class(es):

i. $\mathcal{L}' = \{L_1 \cup \overline{L_1}^A \cup L_2 \cup \overline{L_2}^A\}$ when either $L_1 = L_2$, $L_1 = \overline{L_1}^A$ or $L_2 = \overline{L_2}^A$ ($\overline{L_1}^A = \overline{L_2}^A$ is equivalent to $L_1 = L_2$); in this case the only circuit $L \in \mathcal{L}'$ is its own dual $\overline{L}^{A'} = L$ and $L$ is thus shorted;

ii. $\mathcal{L}' = \{L_1 \cup \overline{L_2}^A,\ \overline{L_1}^A \cup L_2\}$ otherwise; the dualism between the two circuits $\{L, L'\} = \mathcal{L}'$ comes from $e_1 \in L$ and $e_2 \in L'$ becoming duals by construction, while $L \neq L' = \overline{L}^{A'}$ since it would imply the satisfaction of the case (i.) condition.

Since this exhaust the set $\leftrightarrow_{A'} \setminus \leftrightarrow_A$, no other class of $\mathcal{L}_A$ gets merged under $\sim_{A'}$. For the last part, just note that $\mathcal{L}' = \{L_1, \overline{L_1}^A\}$, therefore $\mathcal{L}_{A'} = \{L_1, \overline{L_1}^A\} \cup \mathcal{L}_A \setminus \{L_1, \overline{L_1}^A\} = \mathcal{L}_A$. $\square$

*Proof of Corollary 1.* $\mathcal{L}' = \{L_1, \overline{L_1}^A\}$, therefore $\mathcal{L}_{A'} = \{L_1, \overline{L_1}^A\} \cup \mathcal{L}_A \setminus \{L_1, \overline{L_1}^A\} = \mathcal{L}_A$. $\square$

*Proof of Lemma 2.* Let us consider the linear application

$$\sigma_L^A \colon \langle E_{\overline{\boldsymbol{x}}}^n \rangle \to \mathbb{R}, \qquad \sigma_L^A(\boldsymbol{a}) = \textstyle\sum_{e \in L} a_e - \sum_{e \in \overline{L}^A} a_e,$$

and let $K := \ker \sigma_L^A \cap \{\mathrm{bnd}_e^n \mid e \in L \cup \overline{L}^A\}$. It is easy to see from the definition of $[L]_A$ that $\sigma_L^A([e_1'|e_2']) = 0$ for all $[e_1'|e_2'] \in [L]_{\mathrm{Deg}_{\overline{\boldsymbol{x}}}}$, thus $[L]_{\mathrm{Deg}_{\overline{\boldsymbol{x}}}} \subseteq K$; but since we have that $\dim(K) = 2|L| - 1$ and we already know that $\dim(\mathrm{span}\,[L]_{\mathrm{Deg}_{\overline{\boldsymbol{x}}}}) = 2|L| - 1$, it must be $\mathrm{span}\,[L]_{\mathrm{Deg}_{\overline{\boldsymbol{x}}}} = K$. With regards to the vector $[e_1|e_2]$, remembering that $e_1 \in L$, we have

$$\sigma_L^A([e_1|e_2]) = \begin{cases} 2 & \text{if } e_2 \in L, \\ 1 & \text{if } e_2 \notin L \cup \overline{L}^A, \\ 0 & \text{if } e_2 \in \overline{L}^A. \end{cases}$$

In the last case $[e_1|e_2] \in K = \mathrm{span}\,[L]_{\mathrm{Deg}_{\overline{\boldsymbol{x}}}}$, while in the first two $[e_1|e_2] \notin \mathrm{span}\,[L]_{\mathrm{Deg}_{\overline{\boldsymbol{x}}}}$.

For any $\mathrm{Deg}_{\overline{\boldsymbol{x}}} \subseteq A \subsetneq \mathrm{Deg}_{\overline{\boldsymbol{x}}} \cup \mathrm{Sep}_{\overline{\boldsymbol{x}}}$ it holds $[L]_A = [L]_{\mathrm{Deg}_{\overline{\boldsymbol{x}}}} \cup [L]_{A \setminus \mathrm{Deg}_{\overline{\boldsymbol{x}}}}$, but since

$$\dim(\mathrm{span}\,[L]_A) = \dim(\mathrm{span}\,[L]_{\mathrm{Deg}_{\overline{\boldsymbol{x}}}}) = 2|L| - 1,$$

the linear dependency of $[e_1|e_2]$ from $[L]_A$ boils down to its dependency from $[L]_{\mathrm{Deg}_{\overline{\boldsymbol{x}}}}$. $\square$

*Proof of Lemma 3.* By writing $A = \mathrm{Deg}_{\overline{\boldsymbol{x}}} \cup \{\boldsymbol{a}^1, \ldots, \boldsymbol{a}^m\}$ with $\boldsymbol{a}^k \in \mathrm{Sep}_{\overline{\boldsymbol{x}}}$ it can be seen that $\mathcal{L}_A$ is obtained from $\mathcal{L}_{\mathrm{Deg}_{\overline{\boldsymbol{x}}}}$ by iteratively merging the vectors $(\boldsymbol{a}^1, \ldots, \boldsymbol{a}^m)$ as described in Proposition 1. We thus prove this proposition by induction on the sequence of sets $(A_0, A_1, \ldots, A_m)$ with $A_k := \mathrm{Deg}_{\overline{\boldsymbol{x}}} \cup \bigcup_{i=1}^k \boldsymbol{a}^k,\ k = \{1, ..., m\}$ and $A_0 := \mathrm{Deg}_{\overline{\boldsymbol{x}}}$. The base case $A_0 = \mathrm{Deg}_{\overline{\boldsymbol{x}}}$, where $L \in \mathcal{L}_{A_0}$ is always non-shorted and $\dim(\mathrm{span}\,[L]_{\mathrm{Deg}_{\overline{\boldsymbol{x}}}}) = 2|L| - 1$ has already been proved.

Now for the inductive case. Assume that the thesis holds for $\mathcal{L}_{A_k}$, and given $[e_1|e_2] = \boldsymbol{a}^{k+1}$ take $L_1, L_2 \in \mathcal{L}_{A_k}$ such that $e_1 \in L_1, e_2 \in L_2$. Since $\mathcal{L}_{A_k}$ and

$\mathcal{L}_{A_{k+1}}$ are identical for all elements that are not $L_1, L_2$, their duals and the results of their merging, and since $\operatorname{supp}_{E^n_{\overline{\boldsymbol{x}}}}([e_1|e_2])$ is disjoint from such circuits $L'$ we have that $[L']_{A_{k+1}} = [L']_{A_k}$. It is therefore sufficient to prove the thesis on those $L \in \mathcal{L}'$ resulting from the merging of $L_1, L_2$ with $[e_1|e_2]$. We study separately the different possibilities for $L_1$ and $L_2$.

- ▷ If $L_2 = \overline{L_1}^{A_k}$ then $\mathcal{L}' = \{L_1, \overline{L_1}^{A_k}\}$, but since $e_1 \in L_1, e_2 \in \overline{L_1}^{A_k} = L_2$ Lemma 2 implies $\operatorname{span}[L]_{A_{k+1}} = \operatorname{span}[L]_{A_k}$ for all $L \in \mathcal{L}'$, obtaining this way the thesis.
- ▷ If $L_2 = L_1$ and $L_1$ is not shorted (otherwise also $L_2 = \overline{L_1}^{A_k}$ and it would fall again in the previous case) let $L = L_1 \cup \overline{L_1}^{A_k}$ be the only element of $\mathcal{L}'$, which is shorted in $A_{k+1}$. Then $[e_1|e_2]$ is independent from $[L]_{A_k}$ (from Lemma 2), thus

$$\begin{aligned}\dim(\operatorname{span}[L]_{A_{k+1}}) &= \dim(\operatorname{span}[L_1]_{A_k} \oplus \operatorname{span}\{[e_1|e_2]\}) \\ &= (2|L_1| - 1) + 1 = |L|.\end{aligned}$$

- ▷ If $L_2 \notin \{L_1, \overline{L_1}^{A_k}\}$ and both $L_1$ and $L_2$ are shorted then $L = L_1 \cup L_2$ is again the only element of $\mathcal{L}'$ and is shorted; also since $L_1, L_2$ belong to different pairs they generate disjoint subspaces $\operatorname{span}[L_1]_{A_k}$ and $\operatorname{span}[L_2]_{A_k}$, therefore

$$\dim(\operatorname{span}[L]_{A_{k+1}}) = \dim\left(\operatorname{span}[L_1]_{A_k} \oplus \operatorname{span}[L_2]_{A_k}\right) = |L_1| + |L_2| = |L|.$$

- ▷ If $L_2 \notin \{L_1, \overline{L_1}^{A_k}\}$ and at least one between $L_1$ and $L_2$ is not shorted (say $L_2$), again $\operatorname{span}[L_1]_{A_k}$ and $\operatorname{span}[L_2]_{A_k}$ are disjoint since $L_1, L_2$ belong to different pairs. Consider $\sigma^{A_k}_{L_2}$ as defined in the proof of Lemma 2: while $[e_1|e_2] \notin \ker \sigma^{A_k}_{L_2}$, it is easy to see that $\operatorname{span}[L_1]_{A_k} \oplus \operatorname{span}[L_2]_{A_k} \subseteq \ker \sigma$ since all its generators lie in $\ker \sigma$. Therefore $[e_1|e_2]$ is not contained in $\operatorname{span}[L_1]_{A_k} \oplus \operatorname{span}[L_2]_{A_k}$. Now,
  - ▷ if $L_1$ is shorted then $L = L_1 \cup L_2 \cup \overline{L_2}^{A_k}$ is again the only element of $\mathcal{L}'$ and is shorted, and we have

$$\begin{aligned}\dim(\operatorname{span}[L]_{A_{k+1}}) &= \dim\left(\operatorname{span}[L_1]_{A_k} \oplus \operatorname{span}[L_2]_{A_k} \oplus \operatorname{span}\{[e_1|e_2]\}\right) \\ &= |L_1| + (2|L_2| - 1) + 1 = |L|;\end{aligned}$$

  - ▷ if $L_1$ is also not shorted then $\mathcal{L}' = \left\{L_1 \cup \overline{L_2}^{A}, \overline{L_1}^{A} \cup L_2\right\}$ and both elements of $\mathcal{L}'$ are not shorted. For each $L \in \mathcal{L}'$ we therefore have

$$\begin{aligned}\dim(\operatorname{span}[L]_{A_{k+1}}) &= \dim\left(\operatorname{span}[L_1]_{A_k} \oplus \operatorname{span}[L_2]_{A_k} \oplus \operatorname{span}\{[e_1|e_2]\}\right) \\ &= (2|L_1| - 1) + (2|L_2| - 1) + 1 = 2|L| - 1.\end{aligned}$$

□

*Proof of Proposition 2.* $c$ effectively counts the number of non-shorted circuit pairs in $\mathcal{L}_A$, since a simple case check shows that the number of non-shorted circuits pairs in $\mathcal{L}_A$ decreases by a unit if and only if the condition in point (ii.) is met by $L_1, L_2$. The terminating condition of $c = 0$ is therefore equal to the request that all circuit pairs be shorted. The fact that circuit pairs generate disjoint subspaces which are full rank if and only if the pair is shortened, together with the fact that it is sufficient to check for the rank maximality of the coefficient vectors of active constraints over the set of $1/2$-arcs of $\overline{\boldsymbol{x}}$ proves the correctness of the algorithm. □

# Appendix D Highest integrality gap instances

## $n = 4$

gap: 1.2

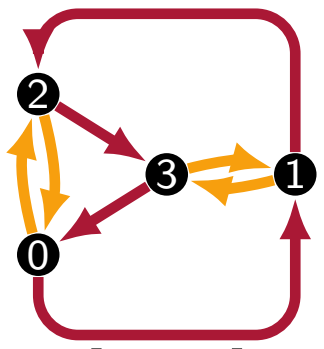

[0123]
[02 | 13]

| ↱ | 0 | 1 | 2 | 3 |
|---|---|---|---|---|
| 0 | | 0 | 1 | 1 |
| 1 | 2 | | 2 | 1 |
| 2 | 0 | 0 | | 0 |
| 3 | 1 | 0 | 2 | |

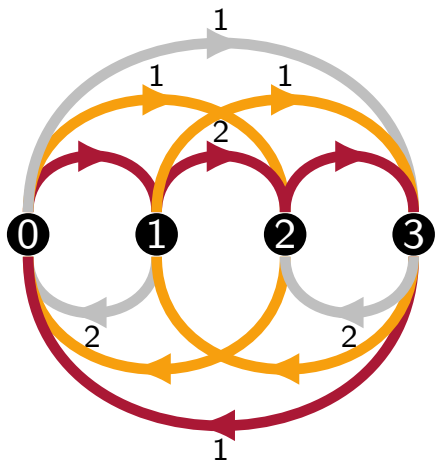

## $n = 5$

gap: 1.25

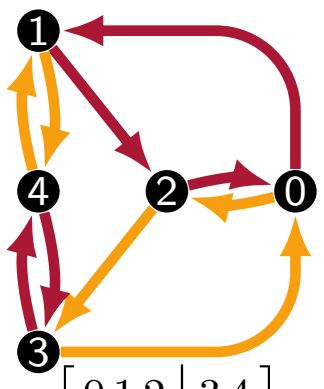

[012 | 34]
[023 | 14]

| ↱ | 0 | 1 | 2 | 3 | 4 |
|---|---|---|---|---|---|
| 0 | | 1 | 2 | 3 | 2 |
| 1 | 2 | | 2 | 2 | 1 |
| 2 | 0 | 1 | | 1 | 1 |
| 3 | 0 | 0 | 2 | | 0 |
| 4 | 1 | 0 | 2 | 1 | |

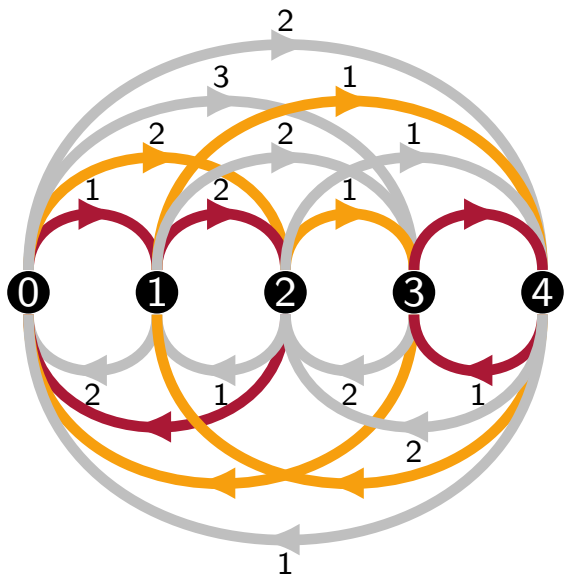

## $n = 6$

gap: 1.33333

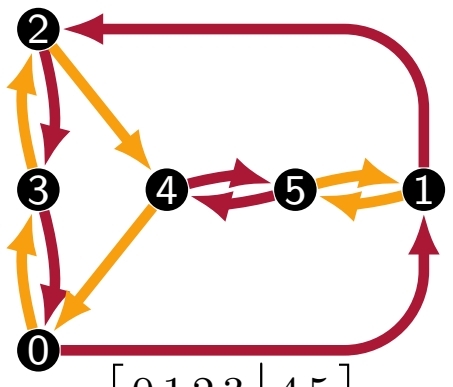

[0123 | 45]
[0324 | 15]

| ↱ | 0 | 1 | 2 | 3 | 4 | 5 |
|---|---|---|---|---|---|---|
| 0 | | 1 | 1 | 1 | 2 | 1 |
| 1 | 1 | | 0 | 1 | 1 | 0 |
| 2 | 1 | 2 | | 1 | 1 | 1 |
| 3 | 0 | 1 | 0 | | 1 | 1 |
| 4 | 0 | 1 | 1 | 1 | | 0 |
| 5 | 1 | 1 | 1 | 2 | 1 | |

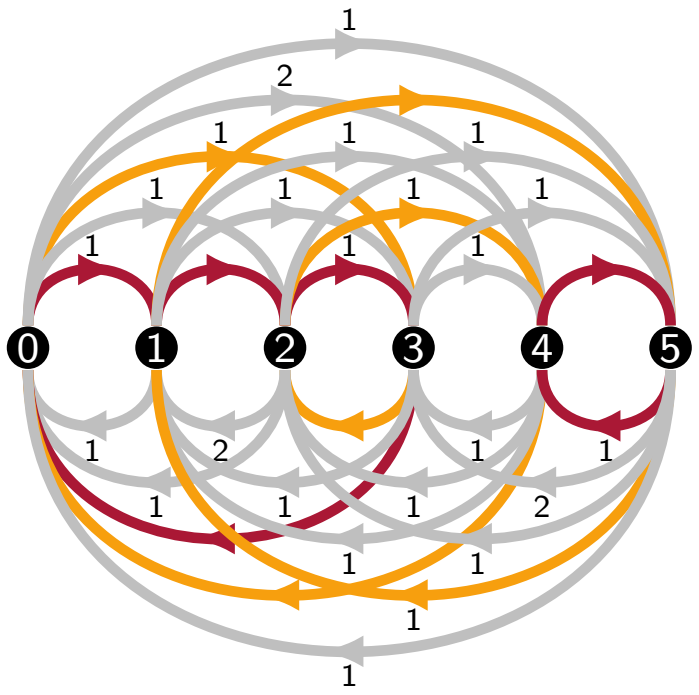

# $n = 7$

gap: 1.33333

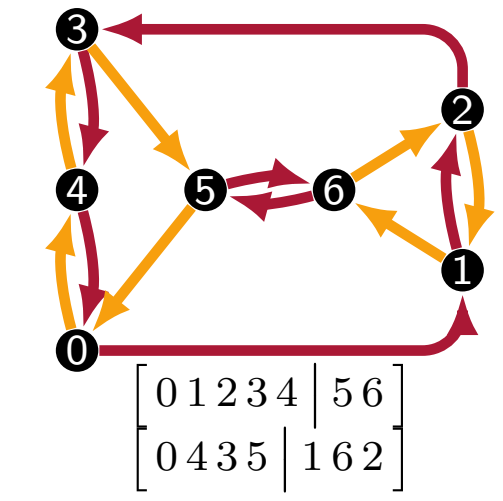

| | 0 | 1 | 2 | 3 | 4 | 5 | 6 |
|---|---|---|---|---|---|---|---|
| 0 | | 0 | 0 | 1 | 0 | 1 | 1 |
| 1 | 2 | | 0 | 1 | 1 | 1 | 1 |
| 2 | 2 | 0 | | 1 | 1 | 1 | 1 |
| 3 | 1 | 1 | 1 | | 0 | 0 | 1 |
| 4 | 1 | 1 | 1 | 1 | | 1 | 2 |
| 5 | 1 | 1 | 1 | 2 | 1 | | 1 |
| 6 | 1 | 0 | 0 | 1 | 1 | 0 | |

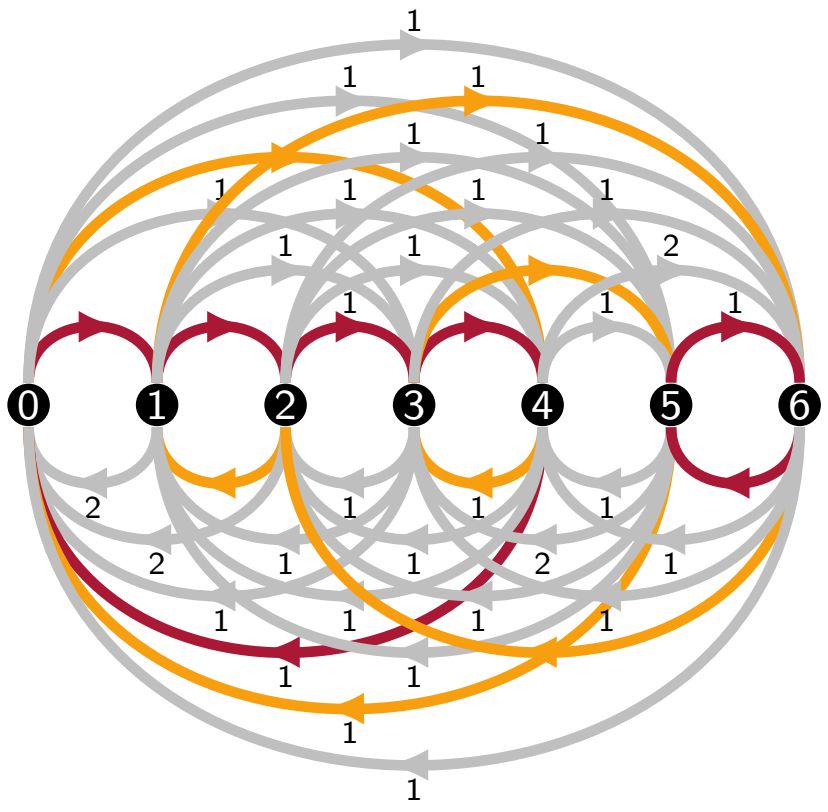

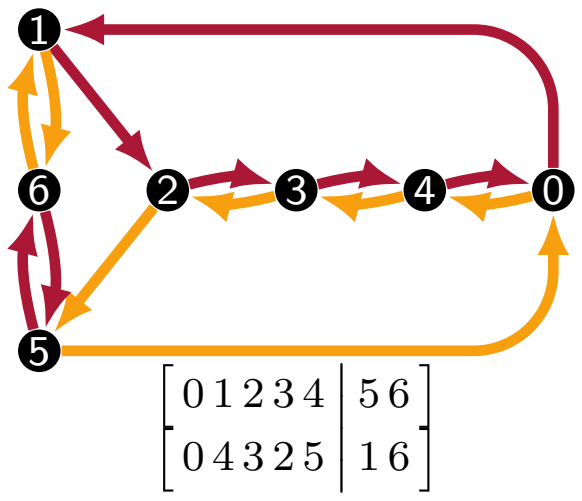

| | 0 | 1 | 2 | 3 | 4 | 5 | 6 |
|---|---|---|---|---|---|---|---|
| 0 | | 1 | 1 | 1 | 0 | 2 | 1 |
| 1 | 1 | | 0 | 0 | 0 | 1 | 0 |
| 2 | 1 | 2 | | 0 | 0 | 1 | 1 |
| 3 | 1 | 2 | 0 | | 0 | 1 | 1 |
| 4 | 1 | 2 | 1 | 1 | | 2 | 2 |
| 5 | 0 | 1 | 1 | 1 | 0 | | 0 |
| 6 | 1 | 1 | 1 | 1 | 1 | 1 | |

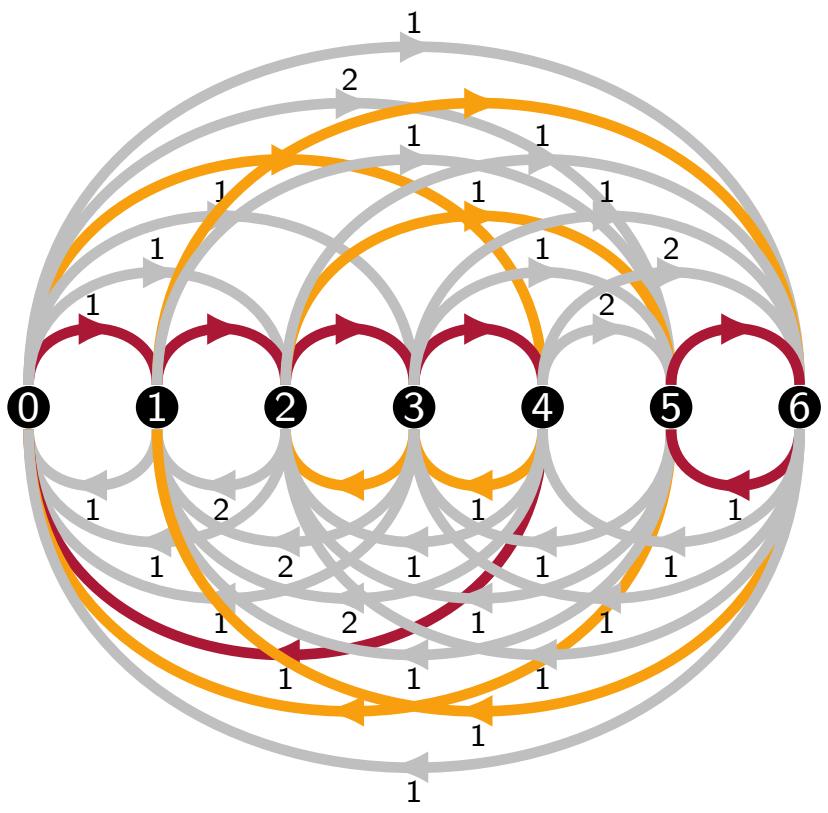

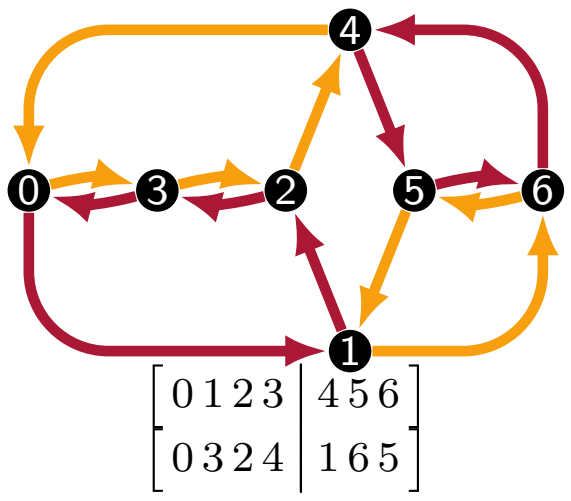

| | 0 | 1 | 2 | 3 | 4 | 5 | 6 |
|---|---|---|---|---|---|---|---|
| 0 | | 1 | 2 | 1 | 2 | 2 | 2 |
| 1 | 1 | | 1 | 1 | 1 | 1 | 1 |
| 2 | 0 | 1 | | 0 | 0 | 1 | 1 |
| 3 | 0 | 1 | 1 | | 1 | 2 | 2 |
| 4 | 0 | 1 | 2 | 1 | | 1 | 1 |
| 5 | 0 | 0 | 1 | 1 | 0 | | 0 |
| 6 | 0 | 0 | 1 | 1 | 0 | 0 | |

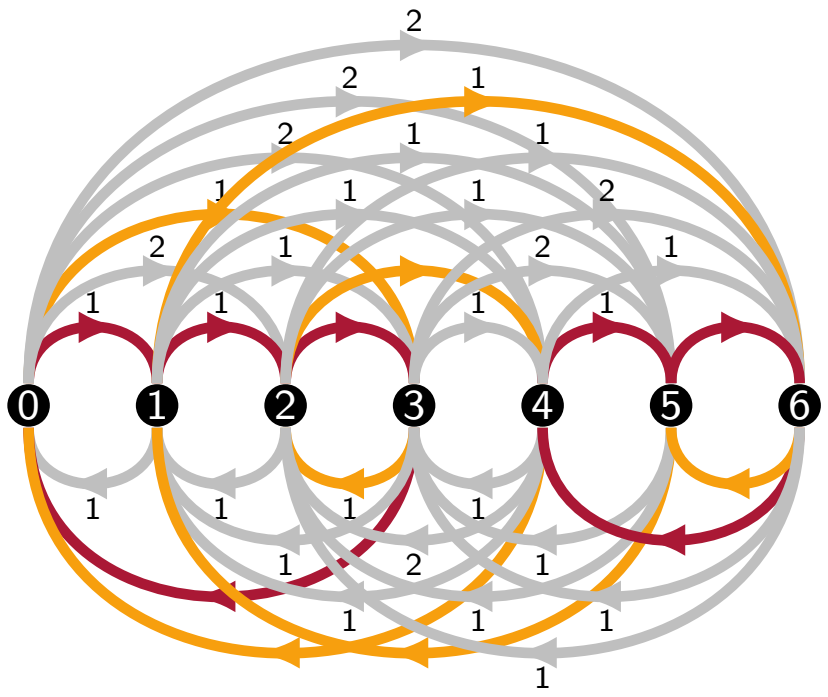

## $n = 8$

`gap:  1.33333`

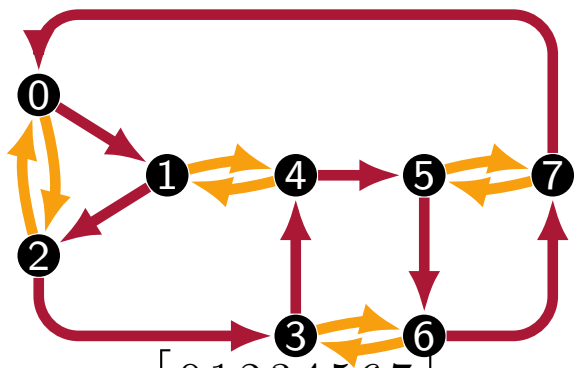

[01234567]
[02|14|36|57]

| ↱ | 0 | 1 | 2 | 3 | 4 | 5 | 6 | 7 |
|---|---|---|---|---|---|---|---|---|
| 0 | | 1 | 0 | 2 | 2 | 2 | 2 | 2 |
| 1 | 0 | | 0 | 2 | 1 | 1 | 2 | 1 |
| 2 | 0 | 1 | | 2 | 2 | 2 | 2 | 2 |
| 3 | 0 | 0 | 0 | | 0 | 0 | 0 | 0 |
| 4 | 0 | 0 | 0 | 2 | | 0 | 1 | 0 |
| 5 | 0 | 1 | 0 | 2 | 2 | | 1 | 0 |
| 6 | 0 | 1 | 0 | 1 | 1 | 0 | | 0 |
| 7 | 0 | 1 | 0 | 2 | 2 | 0 | 1 | |

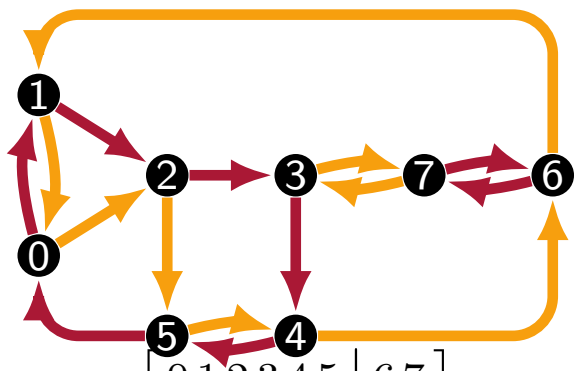

[012345|67]
[025461|37]

| ↱ | 0 | 1 | 2 | 3 | 4 | 5 | 6 | 7 |
|---|---|---|---|---|---|---|---|---|
| 0 | | 0 | 1 | 1 | 1 | 1 | 1 | 1 |
| 1 | 0 | | 1 | 1 | 1 | 1 | 1 | 1 |
| 2 | 0 | 0 | | 0 | 0 | 0 | 1 | 0 |
| 3 | 1 | 1 | 1 | | 0 | 1 | 1 | 0 |
| 4 | 1 | 1 | 2 | 2 | | 1 | 1 | 1 |
| 5 | 0 | 0 | 1 | 1 | 0 | | 1 | 1 |
| 6 | 0 | 0 | 1 | 1 | 1 | 1 | | 0 |
| 7 | 1 | 1 | 2 | 1 | 1 | 2 | 1 | |

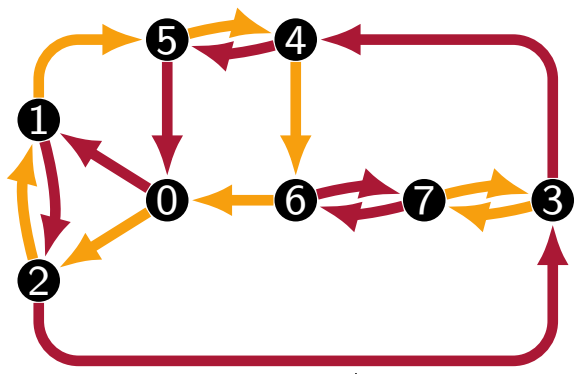

[012345|67]
[021546|37]

| ↱ | 0 | 1 | 2 | 3 | 4 | 5 | 6 | 7 |
|---|---|---|---|---|---|---|---|---|
| 0 | | 0 | 0 | 0 | 0 | 0 | 0 | 0 |
| 1 | 0 | | 0 | 0 | 0 | 0 | 0 | 0 |
| 2 | 0 | 0 | | 0 | 0 | 0 | 0 | 0 |
| 3 | 2 | 2 | 2 | | 0 | 1 | 0 | 0 |
| 4 | 2 | 2 | 2 | 2 | | 1 | 0 | 1 |
| 5 | 1 | 1 | 1 | 1 | 0 | | 0 | 1 |
| 6 | 2 | 2 | 2 | 2 | 2 | 2 | | 1 |
| 7 | 2 | 2 | 2 | 1 | 1 | 2 | 0 | |

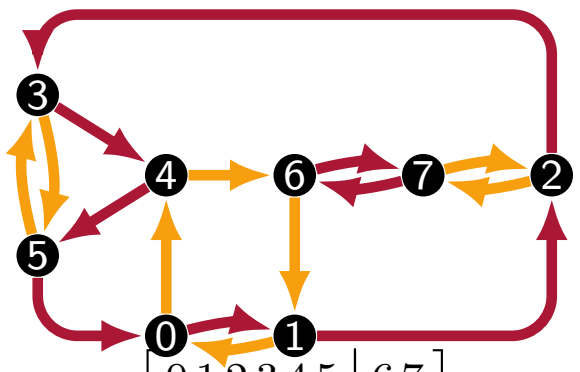

[012345|67]
[0461|27|35]

| ↱ | 0 | 1 | 2 | 3 | 4 | 5 | 6 | 7 |
|---|---|---|---|---|---|---|---|---|
| 0 | | 0 | 2 | 1 | 1 | 1 | 1 | 2 |
| 1 | 1 | | 2 | 2 | 2 | 2 | 2 | 2 |
| 2 | 0 | 0 | | 0 | 0 | 0 | 0 | 0 |
| 3 | 0 | 0 | 2 | | 0 | 0 | 0 | 1 |
| 4 | 0 | 0 | 2 | 0 | | 0 | 0 | 1 |
| 5 | 0 | 0 | 2 | 0 | 0 | | 0 | 1 |
| 6 | 1 | 0 | 2 | 2 | 2 | 2 | | 1 |
| 7 | 1 | 0 | 1 | 1 | 1 | 1 | 0 | |

# $n = 8$

gap: 1.33333

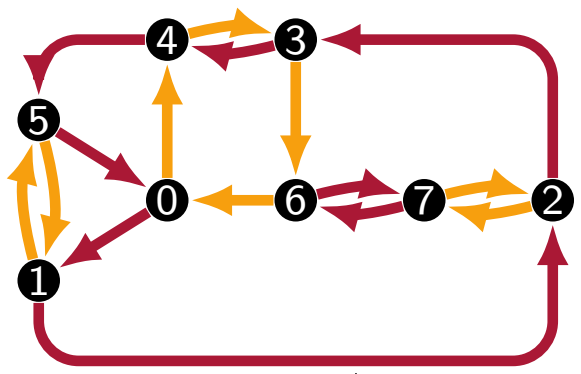

[012345|67]
[0436|15|27]

| ↱ | 0 | 1 | 2 | 3 | 4 | 5 | 6 | 7 |
|---|---|---|---|---|---|---|---|---|
| 0 | | 0 | 0 | 0 | 0 | 0 | 1 | 1 |
| 1 | 0 | | 0 | 0 | 0 | 0 | 1 | 1 |
| 2 | 2 | 2 | | 0 | 1 | 2 | 1 | 1 |
| 3 | 2 | 2 | 2 | | 1 | 2 | 1 | 2 |
| 4 | 1 | 1 | 1 | 0 | | 1 | 1 | 2 |
| 5 | 0 | 0 | 0 | 0 | 0 | | 1 | 1 |
| 6 | 1 | 1 | 1 | 1 | 1 | 1 | | 1 |
| 7 | 1 | 1 | 0 | 0 | 1 | 1 | 0 | |

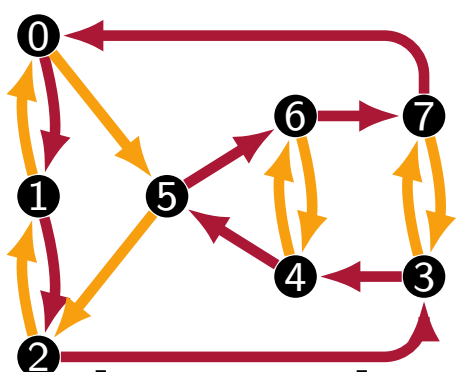

[01234567]
[0521|37|46]

| ↱ | 0 | 1 | 2 | 3 | 4 | 5 | 6 | 7 |
|---|---|---|---|---|---|---|---|---|
| 0 | | 1 | 2 | 2 | 1 | 0 | 1 | 2 |
| 1 | 0 | | 1 | 1 | 1 | 0 | 1 | 1 |
| 2 | 0 | 0 | | 0 | 0 | 0 | 0 | 0 |
| 3 | 0 | 1 | 2 | | 0 | 0 | 0 | 0 |
| 4 | 1 | 2 | 2 | 1 | | 0 | 0 | 1 |
| 5 | 2 | 2 | 2 | 2 | 1 | | 1 | 2 |
| 6 | 1 | 2 | 2 | 1 | 0 | 0 | | 1 |
| 7 | 0 | 1 | 2 | 0 | 0 | 0 | 0 | |

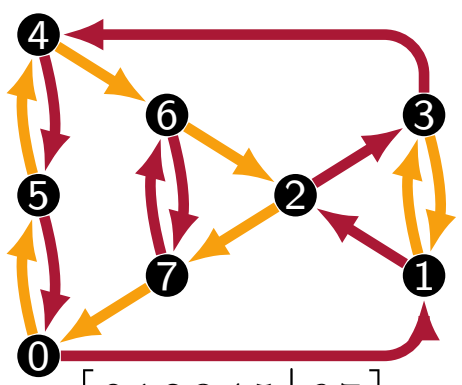

[012345|67]
[054627|13]

| ↱ | 0 | 1 | 2 | 3 | 4 | 5 | 6 | 7 |
|---|---|---|---|---|---|---|---|---|
| 0 | | 0 | 1 | 0 | 1 | 0 | 1 | 1 |
| 1 | 2 | | 1 | 0 | 1 | 1 | 1 | 1 |
| 2 | 1 | 0 | | 0 | 1 | 1 | 0 | 0 |
| 3 | 2 | 0 | 1 | | 1 | 1 | 1 | 1 |
| 4 | 1 | 1 | 1 | 1 | | 0 | 0 | 0 |
| 5 | 1 | 1 | 2 | 1 | 1 | | 1 | 1 |
| 6 | 1 | 1 | 1 | 1 | 2 | 1 | | 0 |
| 7 | 1 | 1 | 1 | 1 | 2 | 1 | 0 | |

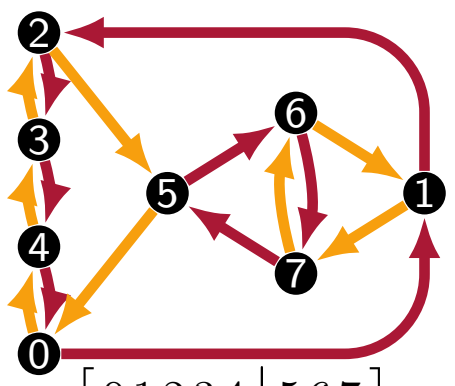

[01234|567]
[04325|176]

| ↱ | 0 | 1 | 2 | 3 | 4 | 5 | 6 | 7 |
|---|---|---|---|---|---|---|---|---|
| 0 | | 2 | 2 | 2 | 1 | 2 | 2 | 2 |
| 1 | 0 | | 0 | 0 | 0 | 0 | 0 | 0 |
| 2 | 0 | 2 | | 0 | 0 | 0 | 1 | 1 |
| 3 | 0 | 2 | 0 | | 0 | 0 | 1 | 1 |
| 4 | 0 | 2 | 1 | 1 | | 1 | 2 | 2 |
| 5 | 0 | 2 | 2 | 2 | 1 | | 1 | 1 |
| 6 | 0 | 1 | 1 | 1 | 1 | 0 | | 0 |
| 7 | 0 | 1 | 1 | 1 | 1 | 0 | 0 | |

## $n = 8$

gap:  1.33333

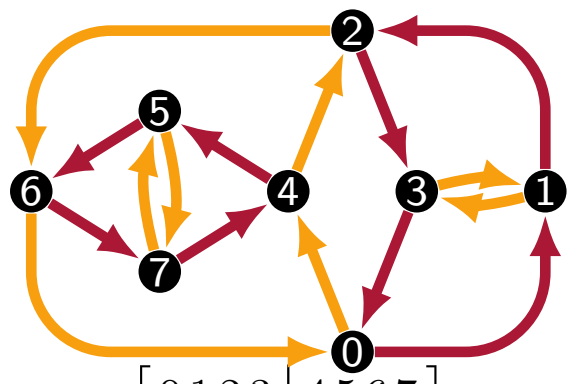

[0123 | 4567]
[0426 | 13 | 57]

| ↱ | 0 | 1 | 2 | 3 | 4 | 5 | 6 | 7 |
|---|---|---|---|---|---|---|---|---|
| 0 | | 0 | 0 | 0 | 0 | 1 | 2 | 1 |
| 1 | 1 | | 0 | 0 | 1 | 2 | 2 | 2 |
| 2 | 2 | 1 | | 1 | 2 | 2 | 2 | 2 |
| 3 | 1 | 0 | 0 | | 1 | 2 | 2 | 2 |
| 4 | 2 | 1 | 0 | 1 | | 1 | 2 | 1 |
| 5 | 1 | 1 | 0 | 1 | 0 | | 1 | 0 |
| 6 | 0 | 0 | 0 | 0 | 0 | 0 | | 0 |
| 7 | 1 | 1 | 0 | 1 | 0 | 0 | 1 | |

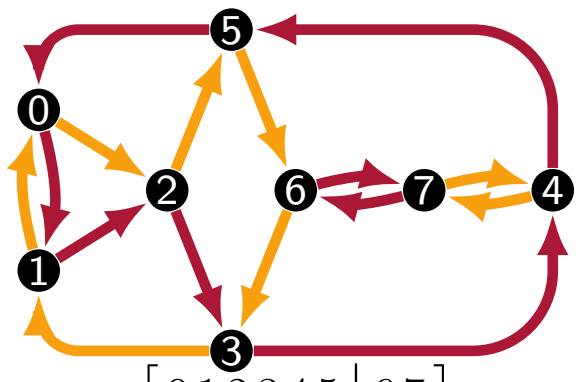

[012345 | 67]
[025631 | 47]

| ↱ | 0 | 1 | 2 | 3 | 4 | 5 | 6 | 7 |
|---|---|---|---|---|---|---|---|---|
| 0 | | 0 | 1 | 1 | 1 | 1 | 1 | 1 |
| 1 | 0 | | 1 | 1 | 1 | 1 | 1 | 1 |
| 2 | 0 | 0 | | 0 | 0 | 0 | 0 | 1 |
| 3 | 0 | 0 | 1 | | 0 | 1 | 1 | 1 |
| 4 | 1 | 1 | 2 | 2 | | 1 | 1 | 1 |
| 5 | 0 | 0 | 1 | 1 | 1 | | 0 | 1 |
| 6 | 1 | 1 | 2 | 1 | 1 | 2 | | 1 |
| 7 | 1 | 1 | 1 | 1 | 0 | 1 | 0 | |

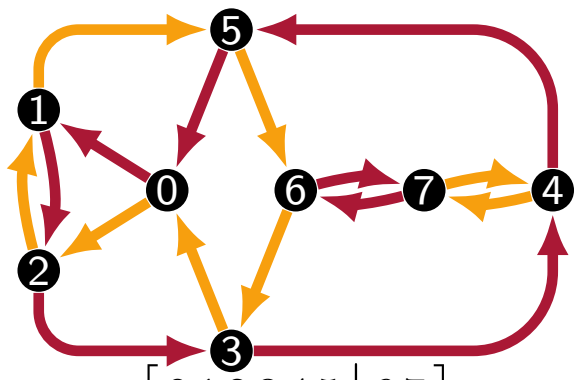

[012345 | 67]
[021563 | 47]

| ↱ | 0 | 1 | 2 | 3 | 4 | 5 | 6 | 7 |
|---|---|---|---|---|---|---|---|---|
| 0 | | 0 | 0 | 0 | 2 | 1 | 1 | 1 |
| 1 | 1 | | 0 | 0 | 2 | 1 | 1 | 2 |
| 2 | 1 | 0 | | 0 | 2 | 1 | 1 | 2 |
| 3 | 1 | 1 | 1 | | 2 | 2 | 2 | 2 |
| 4 | 0 | 0 | 0 | 0 | | 0 | 0 | 0 |
| 5 | 0 | 0 | 0 | 0 | 2 | | 0 | 1 |
| 6 | 1 | 1 | 1 | 0 | 2 | 2 | | 1 |
| 7 | 1 | 0 | 0 | 0 | 1 | 1 | 0 | |

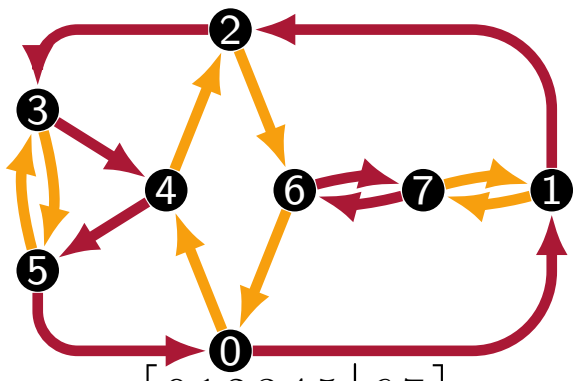

[012345 | 67]
[0426 | 17 | 35]

| ↱ | 0 | 1 | 2 | 3 | 4 | 5 | 6 | 7 |
|---|---|---|---|---|---|---|---|---|
| 0 | | 0 | 1 | 1 | 1 | 1 | 1 | 1 |
| 1 | 2 | | 1 | 2 | 2 | 2 | 1 | 1 |
| 2 | 1 | 1 | | 1 | 1 | 1 | 0 | 1 |
| 3 | 0 | 0 | 0 | | 0 | 0 | 0 | 1 |
| 4 | 0 | 0 | 0 | 0 | | 0 | 0 | 1 |
| 5 | 0 | 0 | 0 | 0 | 0 | | 0 | 1 |
| 6 | 1 | 1 | 2 | 2 | 2 | 2 | | 1 |
| 7 | 1 | 0 | 1 | 2 | 2 | 2 | 0 | |

## $n = 8$

gap: 1.33333

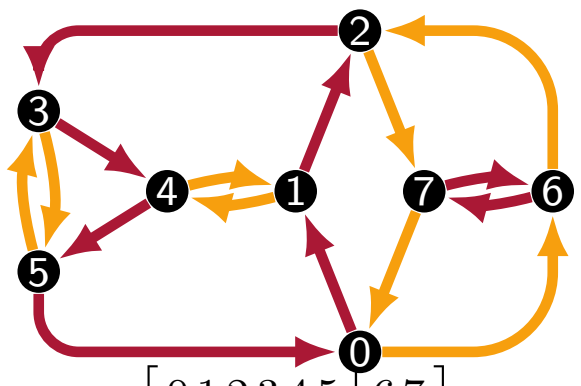

[012345|67]
[0627|14|35]

| | 0 | 1 | 2 | 3 | 4 | 5 | 6 | 7 |
|---|---|---|---|---|---|---|---|---|
| 0 | | 2 | 2 | 2 | 2 | 2 | 1 | 1 |
| 1 | 0 | | 0 | 0 | 0 | 0 | 0 | 0 |
| 2 | 0 | 2 | | 0 | 1 | 0 | 0 | 0 |
| 3 | 0 | 2 | 2 | | 1 | 0 | 1 | 1 |
| 4 | 0 | 1 | 1 | 0 | | 0 | 1 | 1 |
| 5 | 0 | 2 | 2 | 0 | 1 | | 1 | 1 |
| 6 | 0 | 2 | 1 | 1 | 2 | 1 | | 0 |
| 7 | 0 | 2 | 1 | 1 | 2 | 1 | 0 | |

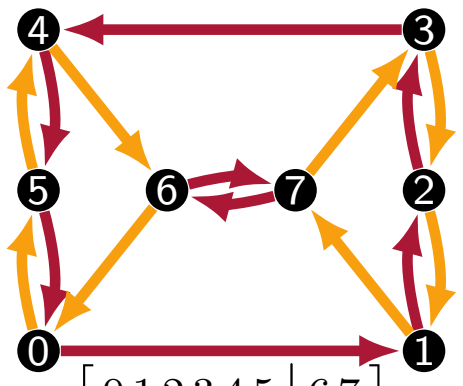

[012345|67]
[0546|1732]

| | 0 | 1 | 2 | 3 | 4 | 5 | 6 | 7 |
|---|---|---|---|---|---|---|---|---|
| 0 | | 1 | 1 | 1 | 2 | 1 | 2 | 1 |
| 1 | 1 | | 0 | 0 | 1 | 1 | 1 | 0 |
| 2 | 1 | 0 | | 0 | 1 | 1 | 1 | 0 |
| 3 | 1 | 0 | 0 | | 1 | 1 | 1 | 0 |
| 4 | 0 | 1 | 1 | 1 | | 0 | 0 | 0 |
| 5 | 0 | 1 | 1 | 1 | 1 | | 1 | 1 |
| 6 | 0 | 1 | 1 | 1 | 2 | 1 | | 0 |
| 7 | 1 | 1 | 1 | 1 | 2 | 2 | 1 | |

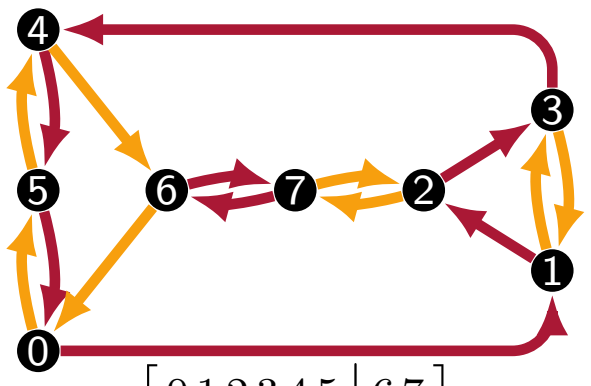

[012345|67]
[0546|13|27]

| | 0 | 1 | 2 | 3 | 4 | 5 | 6 | 7 |
|---|---|---|---|---|---|---|---|---|
| 0 | | 0 | 0 | 0 | 0 | 0 | 0 | 0 |
| 1 | 2 | | 0 | 0 | 0 | 1 | 0 | 0 |
| 2 | 2 | 1 | | 1 | 1 | 2 | 0 | 0 |
| 3 | 2 | 0 | 0 | | 0 | 1 | 0 | 0 |
| 4 | 2 | 2 | 1 | 2 | | 1 | 0 | 1 |
| 5 | 1 | 1 | 1 | 1 | 0 | | 0 | 1 |
| 6 | 2 | 2 | 1 | 2 | 2 | 2 | | 1 |
| 7 | 2 | 1 | 0 | 1 | 1 | 2 | 0 | |

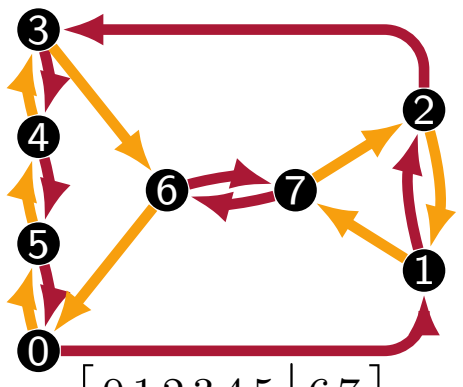

[012345|67]
[05436|172]

| | 0 | 1 | 2 | 3 | 4 | 5 | 6 | 7 |
|---|---|---|---|---|---|---|---|---|
| 0 | | 1 | 1 | 1 | 0 | 0 | 2 | 2 |
| 1 | 1 | | 0 | 0 | 0 | 0 | 1 | 1 |
| 2 | 1 | 0 | | 0 | 0 | 0 | 1 | 1 |
| 3 | 1 | 2 | 2 | | 0 | 0 | 1 | 2 |
| 4 | 1 | 2 | 2 | 1 | | 0 | 2 | 3 |
| 5 | 1 | 2 | 2 | 1 | 0 | | 2 | 3 |
| 6 | 0 | 1 | 1 | 1 | 0 | 0 | | 1 |
| 7 | 0 | 0 | 0 | 0 | 0 | 0 | 0 | |

## $n = 8$

gap: 1.33333

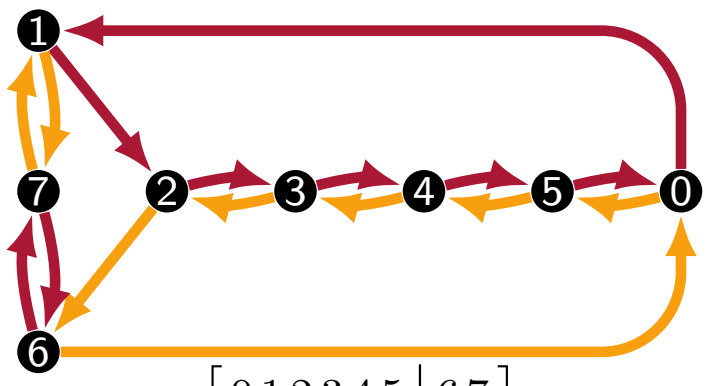

[012345 | 67]
[054326 | 17]

| | 0 | 1 | 2 | 3 | 4 | 5 | 6 | 7 |
|---|---|---|---|---|---|---|---|---|
| 0 | | 0 | 1 | 1 | 1 | 1 | 2 | 1 |
| 1 | 2 | | 1 | 2 | 2 | 2 | 2 | 1 |
| 2 | 1 | 1 | | 1 | 1 | 1 | 1 | 1 |
| 3 | 0 | 0 | 0 | | 0 | 0 | 1 | 1 |
| 4 | 0 | 0 | 0 | 0 | | 0 | 1 | 1 |
| 5 | 0 | 0 | 0 | 0 | 0 | | 1 | 1 |
| 6 | 0 | 0 | 1 | 1 | 1 | 1 | | 0 |
| 7 | 1 | 0 | 1 | 2 | 2 | 2 | 1 | |

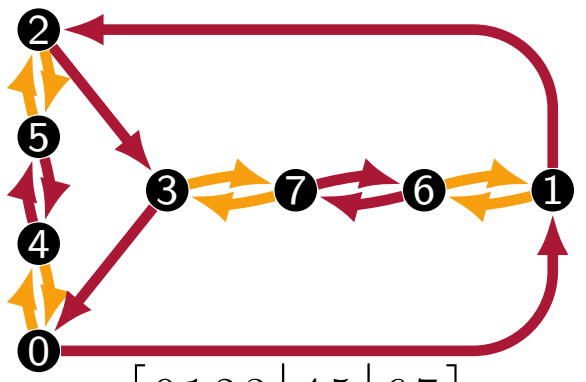

[0123 | 45 | 67]
[04 | 16 | 25 | 37]

| | 0 | 1 | 2 | 3 | 4 | 5 | 6 | 7 |
|---|---|---|---|---|---|---|---|---|
| 0 | | 0 | 0 | 1 | 0 | 0 | 0 | 0 |
| 1 | 2 | | 0 | 1 | 1 | 1 | 0 | 0 |
| 2 | 2 | 2 | | 1 | 1 | 1 | 2 | 1 |
| 3 | 1 | 1 | 1 | | 1 | 1 | 1 | 0 |
| 4 | 1 | 1 | 0 | 1 | | 0 | 1 | 1 |
| 5 | 1 | 1 | 0 | 1 | 0 | | 1 | 1 |
| 6 | 2 | 0 | 0 | 1 | 1 | 1 | | 0 |
| 7 | 2 | 1 | 1 | 1 | 2 | 2 | 1 | |

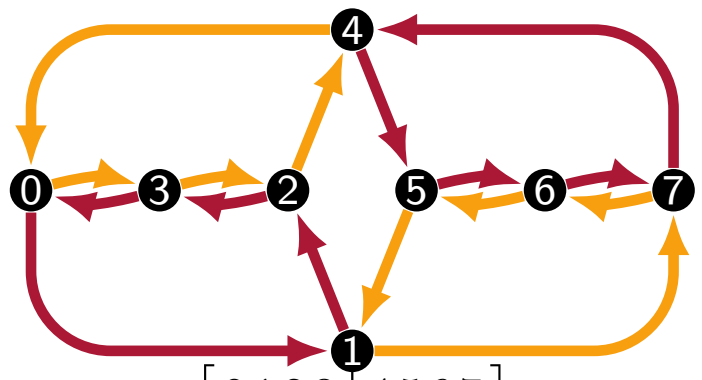

[0123 | 4567]
[0324 | 1765]

| | 0 | 1 | 2 | 3 | 4 | 5 | 6 | 7 |
|---|---|---|---|---|---|---|---|---|
| 0 | | 0 | 1 | 1 | 1 | 1 | 1 | 1 |
| 1 | 2 | | 1 | 2 | 1 | 1 | 1 | 1 |
| 2 | 1 | 1 | | 1 | 0 | 1 | 1 | 1 |
| 3 | 0 | 0 | 0 | | 0 | 1 | 1 | 1 |
| 4 | 1 | 1 | 2 | 2 | | 1 | 1 | 1 |
| 5 | 1 | 0 | 1 | 2 | 0 | | 0 | 0 |
| 6 | 1 | 0 | 1 | 2 | 0 | 0 | | 0 |
| 7 | 1 | 0 | 1 | 2 | 0 | 0 | 0 | |

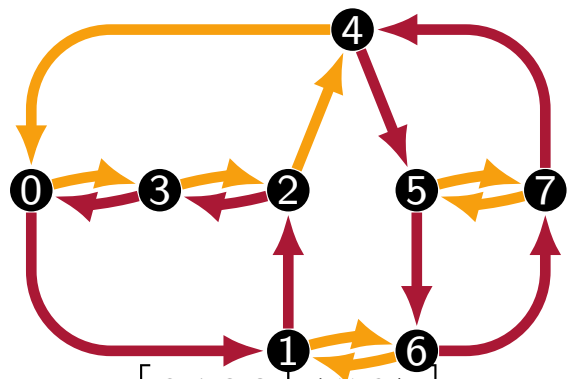

[0123 | 4567]
[0324 | 16 | 57]

| | 0 | 1 | 2 | 3 | 4 | 5 | 6 | 7 |
|---|---|---|---|---|---|---|---|---|
| 0 | | 4 | 4 | 2 | 4 | 5 | 4 | 6 |
| 1 | 0 | | 0 | 0 | 0 | 1 | 0 | 2 |
| 2 | 0 | 4 | | 0 | 0 | 2 | 3 | 4 |
| 3 | 0 | 4 | 2 | | 2 | 4 | 4 | 6 |
| 4 | 0 | 4 | 4 | 2 | | 2 | 3 | 4 |
| 5 | 2 | 4 | 4 | 4 | 2 | | 2 | 2 |
| 6 | 1 | 2 | 2 | 2 | 2 | 2 | | 2 |
| 7 | 0 | 3 | 2 | 2 | 0 | 0 | 2 | |

## $n = 9$

gap: 1.375

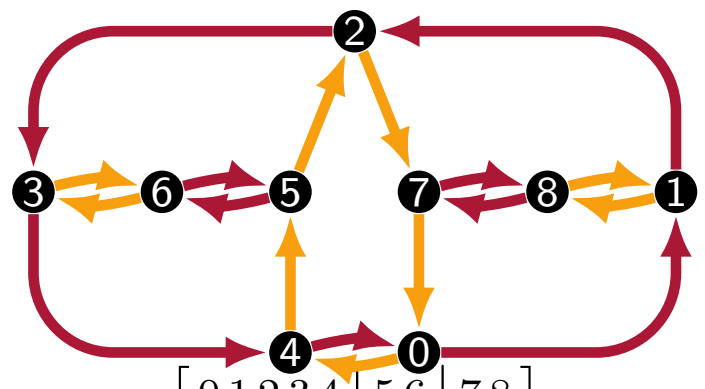

[ 01234 | 56 | 78 ]
[ 04527 | 18 | 36 ]

| ↱ | 0 | 1 | 2 | 3 | 4 | 5 | 6 | 7 | 8 |
|---|---|---|---|---|---|---|---|---|---|
| 0 | | 0 | 4 | 5 | 2 | 3 | 4 | 4 | 2 |
| 1 | 4 | | 4 | 5 | 4 | 5 | 5 | 4 | 2 |
| 2 | 0 | 0 | | 1 | 1 | 1 | 1 | 0 | 0 |
| 3 | 0 | 0 | 2 | | 0 | 0 | 0 | 2 | 2 |
| 4 | 0 | 0 | 3 | 3 | | 1 | 2 | 2 | 2 |
| 5 | 2 | 2 | 2 | 2 | 2 | | 1 | 2 | 2 |
| 6 | 1 | 1 | 2 | 1 | 1 | 0 | | 2 | 2 |
| 7 | 0 | 0 | 4 | 5 | 2 | 3 | 4 | | 0 |
| 8 | 2 | 0 | 4 | 5 | 4 | 5 | 5 | 2 | |

## $n = 10$

gap: 1.4

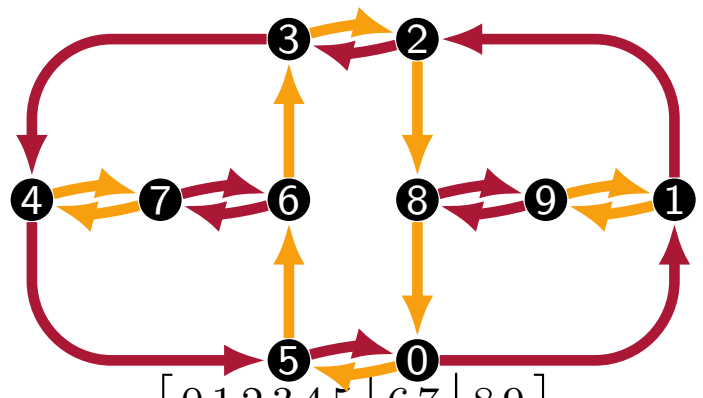

[ 012345 | 67 | 89 ]
[ 056328 | 19 | 47 ]

| ↱ | 0 | 1 | 2 | 3 | 4 | 5 | 6 | 7 | 8 | 9 |
|---|---|---|---|---|---|---|---|---|---|---|
| 0 | | 0 | 1 | 1 | 3 | 1 | 1 | 2 | 1 | 0 |
| 1 | 2 | | 1 | 1 | 3 | 3 | 2 | 3 | 1 | 0 |
| 2 | 1 | 1 | | 0 | 2 | 2 | 2 | 2 | 0 | 0 |
| 3 | 2 | 2 | 1 | | 2 | 2 | 2 | 2 | 1 | 1 |
| 4 | 0 | 0 | 1 | 0 | | 0 | 0 | 0 | 0 | 0 |
| 5 | 0 | 0 | 1 | 0 | 2 | | 0 | 1 | 1 | 0 |
| 6 | 2 | 1 | 1 | 0 | 2 | 2 | | 1 | 1 | 1 |
| 7 | 1 | 1 | 1 | 0 | 1 | 1 | 0 | | 1 | 1 |
| 8 | 1 | 1 | 2 | 2 | 3 | 2 | 2 | 3 | | 0 |
| 9 | 2 | 1 | 2 | 2 | 4 | 3 | 3 | 4 | 1 | |

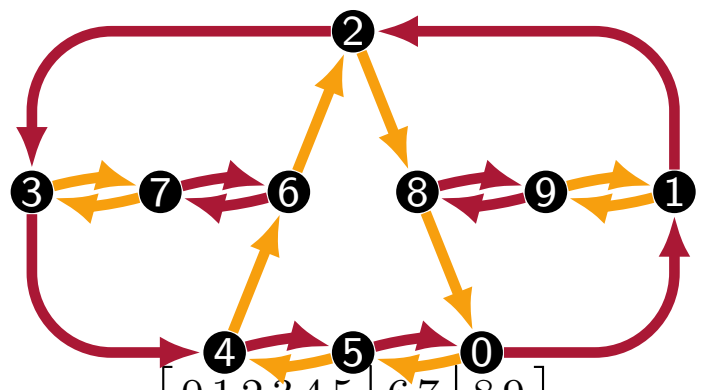

[ 012345 | 67 | 89 ]
[ 054628 | 19 | 37 ]

| ↱ | 0 | 1 | 2 | 3 | 4 | 5 | 6 | 7 | 8 | 9 |
|---|---|---|---|---|---|---|---|---|---|---|
| 0 | | 0 | 1 | 2 | 1 | 0 | 1 | 2 | 2 | 1 |
| 1 | 2 | | 1 | 2 | 2 | 1 | 2 | 2 | 2 | 1 |
| 2 | 1 | 1 | | 1 | 1 | 0 | 1 | 1 | 1 | 1 |
| 3 | 1 | 1 | 1 | | 0 | 0 | 0 | 0 | 2 | 2 |
| 4 | 1 | 1 | 1 | 2 | | 0 | 0 | 1 | 2 | 2 |
| 5 | 1 | 1 | 1 | 2 | 1 | | 1 | 2 | 2 | 2 |
| 6 | 2 | 2 | 1 | 2 | 2 | 1 | | 1 | 2 | 2 |
| 7 | 2 | 2 | 1 | 1 | 1 | 1 | 0 | | 2 | 2 |
| 8 | 0 | 0 | 1 | 2 | 1 | 0 | 1 | 2 | | 0 |
| 9 | 1 | 0 | 1 | 2 | 2 | 1 | 2 | 2 | 1 | |

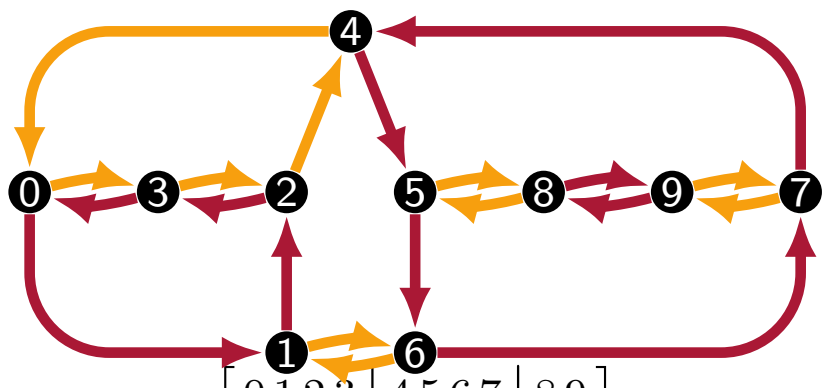

[ 0123 | 4567 | 89 ]
[ 0324 | 16 | 58 | 79 ]

| ↱ | 0 | 1 | 2 | 3 | 4 | 5 | 6 | 7 | 8 | 9 |
|---|---|---|---|---|---|---|---|---|---|---|
| 0 | | 1 | 0 | 0 | 3 | 3 | 5 | 5 | 5 | 5 |
| 1 | 6 | | 0 | 3 | 3 | 2 | 4 | 4 | 4 | 4 |
| 2 | 6 | 7 | | 3 | 3 | 3 | 9 | 8 | 5 | 7 |
| 3 | 3 | 4 | 0 | | 3 | 3 | 8 | 8 | 5 | 7 |
| 4 | 4 | 5 | 4 | 4 | | 0 | 6 | 6 | 2 | 4 |
| 5 | 9 | 6 | 6 | 9 | 6 | | 6 | 6 | 2 | 4 |
| 6 | 4 | 0 | 0 | 3 | 0 | 0 | | 0 | 0 | 0 |
| 7 | 4 | 4 | 4 | 4 | 0 | 0 | 6 | | 0 | 0 |
| 8 | 8 | 6 | 6 | 8 | 4 | 0 | 6 | 4 | | 2 |
| 9 | 6 | 6 | 6 | 6 | 2 | 0 | 6 | 2 | 0 | |

## $n = 11$

gap: 1.42857

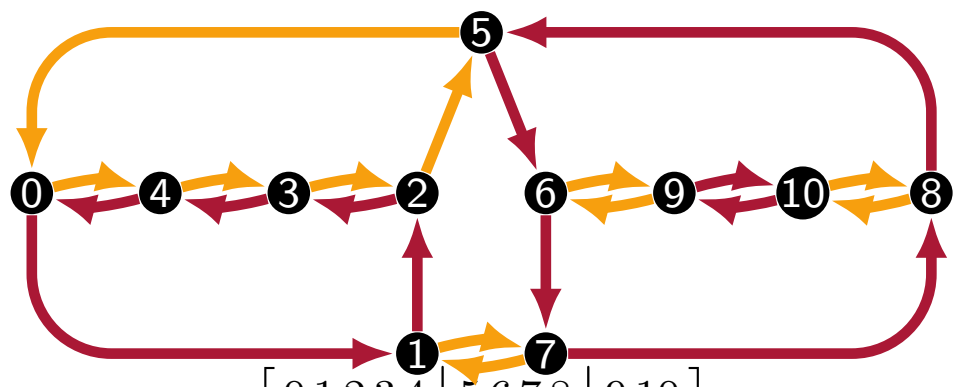

$\left[0\,1\,2\,3\,4 \,\middle|\, 5\,6\,7\,8 \,\middle|\, 9\,10\right]$

$\left[0\,4\,3\,2\,5 \,\middle|\, 1\,7 \,\middle|\, 6\,9 \,\middle|\, 8\,10\right]$

| ↱ | 0 | 1 | 2 | 3 | 4 | 5 | 6 | 7 | 8 | 9 | 10 |
|---|---|---|---|---|---|---|---|---|---|---|---|
| 0 | | 0 | 2 | 1 | 1 | 3 | 3 | 1 | 3 | 4 | 3 |
| 1 | 3 | | 2 | 2 | 3 | 3 | 3 | 1 | 3 | 4 | 3 |
| 2 | 1 | 1 | | 0 | 1 | 1 | 1 | 1 | 3 | 2 | 2 |
| 3 | 1 | 1 | 1 | | 1 | 2 | 2 | 2 | 4 | 3 | 3 |
| 4 | 0 | 0 | 1 | 0 | | 2 | 2 | 1 | 3 | 3 | 3 |
| 5 | 0 | 0 | 2 | 1 | 1 | | 0 | 0 | 2 | 1 | 1 |
| 6 | 3 | 1 | 3 | 3 | 4 | 3 | | 0 | 2 | 1 | 1 |
| 7 | 3 | 1 | 3 | 3 | 4 | 3 | 3 | | 2 | 3 | 2 |
| 8 | 1 | 1 | 3 | 2 | 2 | 1 | 1 | 1 | | 1 | 0 |
| 9 | 2 | 1 | 3 | 3 | 3 | 2 | 0 | 0 | 1 | | 0 |
| 10 | 2 | 2 | 4 | 3 | 3 | 2 | 1 | 1 | 1 | 1 | |

## $n = 12$

gap: 1.43589

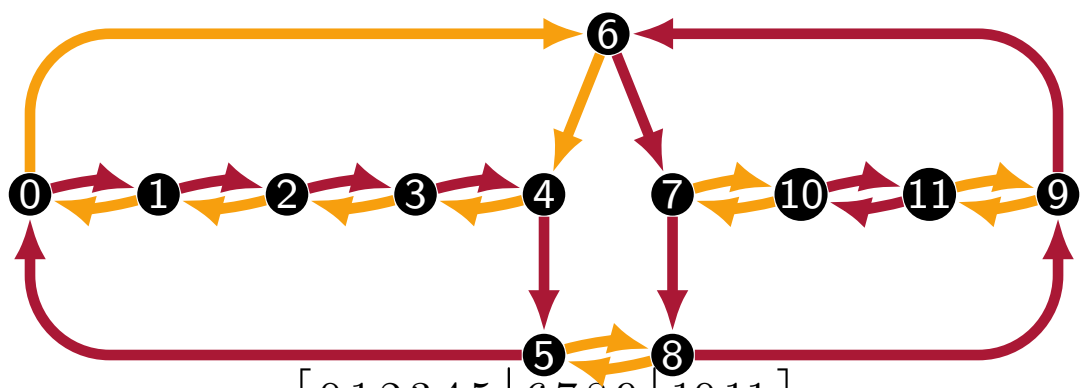

$\left[0\,1\,2\,3\,4\,5 \,\middle|\, 6\,7\,8\,9 \,\middle|\, 10\,11\right]$

$\left[0\,6\,4\,3\,2\,1 \,\middle|\, 5\,8 \,\middle|\, 7\,10 \,\middle|\, 9\,11\right]$

| ↱ | 0 | 1 | 2 | 3 | 4 | 5 | 6 | 7 | 8 | 9 | 10 | 11 |
|---|---|---|---|---|---|---|---|---|---|---|---|---|
| 0 | | 0 | 2 | 4 | 4 | 8 | 4 | 4 | 7 | 7 | 4 | 7 |
| 1 | 2 | | 2 | 4 | 4 | 8 | 6 | 6 | 9 | 9 | 6 | 9 |
| 2 | 2 | 0 | | 2 | 2 | 6 | 6 | 6 | 7 | 7 | 6 | 8 |
| 3 | 2 | 0 | 0 | | 0 | 4 | 6 | 6 | 5 | 5 | 6 | 8 |
| 4 | 4 | 2 | 2 | 2 | | 4 | 8 | 8 | 5 | 5 | 8 | 8 |
| 5 | 0 | 0 | 2 | 4 | 4 | | 4 | 4 | 1 | 1 | 4 | 4 |
| 6 | 4 | 2 | 2 | 2 | 0 | 4 | | 0 | 3 | 3 | 0 | 3 |
| 7 | 7 | 7 | 9 | 11 | 9 | 7 | 9 | | 3 | 3 | 0 | 3 |
| 8 | 4 | 4 | 6 | 8 | 6 | 4 | 6 | 6 | | 0 | 3 | 3 |
| 9 | 10 | 8 | 8 | 8 | 6 | 10 | 6 | 6 | 9 | | 3 | 3 |
| 10 | 10 | 10 | 10 | 11 | 9 | 10 | 9 | 3 | 6 | 3 | | 3 |
| 11 | 10 | 8 | 8 | 8 | 6 | 10 | 6 | 3 | 6 | 0 | 0 | |